\documentclass[12pt,a4paper,oneside]{article}
\usepackage{graphicx,epsfig}
\usepackage[all]{xy}
\usepackage{amssymb,amsmath,amsfonts,amsthm,latexsym}
\usepackage[T1]{fontenc}

\newcommand{\ba}{\begin{array}}
\newcommand{\bea}{\begin{eqnarray}}
\newcommand{\bean}{\begin{eqnarray*}}
\newcommand{\be}{\begin{equation}}
\newcommand{\bl}{\ba{llll}}

\newcommand{\CaSz}{Cauchy-Schwarz inequality~}

\newcommand{\D}{\mathcal{D}}
\def \ddt(#1){\stl{\frac{\partial #1}{\partial t}}}
\newcommand{\Dir}{Dirichlet~}
\newcommand{\dv}{\nabla\cdot}
\newcommand{\dx}{~d\bold{x}}
\newcommand{\dt}{\delta t}

\newcommand{\E}{\mat{E}}
\newcommand{\ea}{\end{array}}
\newcommand{\eea}{\end{eqnarray}}
\newcommand{\eean}{\end{eqnarray*}}
\newcommand{\ee}{\end{equation}}
\newcommand{\el}{\ea}

\newcommand{\grad}{\nabla}

\newcommand{\half}{\frac{1}{2}}
\newcommand{\hD}{h_{\D}}
\newcommand{\hyp}{\mathcal{H}}

\newcommand{\Itn}{\int^{t_n}_{t_{n-1}}}
\newcommand{\IK}{\int_{K}}
\newcommand{\ItnK}{\Itn \IK}
\newcommand{\IT}{\int^{T}_{0}}
\newcommand{\IO}{\int_{\Om}}
\newcommand{\ITO}{\IT \IO}

\newcommand{\intpos}{\mathbb{N}\setminus\{0\}}

\newcommand{\Lam}{\bold{\Lambda}}
\newcommand{\lc}{\lceil}
\newcommand{\LO}{L^2(\Om)}
\newcommand{\lqn}{\lefteqn}
\newcommand{\LQT}{L^2(Q_T)}
\newcommand{\LinfO}{L^{\infty}(0,T;\LO)}
\newcommand{\Lip}{Lipschitz~}

\newcommand{\M}{\mat{M}}
\newcommand{\mat}{\mathcal}

\newcommand{\n}{\bold{n}}
\newcommand{\nKs}{\n_{K,\s}}
\def \norm(#1){\|#1\|}

\newcommand{\Om}{\Omega}

\newcommand{\ph}{\varphi}
\newcommand{\Pb}{(\mat{P})}

\newcommand{\rc}{\rceil}
\newcommand{\real}{\mathbb{R}}
\newcommand{\rhs}{right-hand side~}

\newcommand{\s}{\sigma}
\newcommand{\SK}{\sum_{K\in\M}}
\newcommand{\Sn}{\sum^{N}_{n=1}}
\newcommand{\SnK}{\Sn\SK}
\newcommand{\SnKs}{\Sn\SK\Ss}
\newcommand{\Ss}{\sum_{\s\in\E_K}}
\newcommand{\stl}{\displaystyle}

\newcommand{\ta}{\text{~and~}}
\newcommand{\tf}{\text{~for~}}
\newcommand{\tfa}{\text{~for all~}}

\newcommand{\uhk}{u_{\D,\dt}}
\newcommand{\us}{u_{\s}}
\newcommand{\uK}{u_{K}}
\newcommand{\ubar}{\overline{u^n_{K,\s}}}

\newcommand{\V}{\bold{V}}
\newcommand{\Vmax}{\|\V\|_{L^{\infty}(\Omega)}}

\newcommand{\w}{\bold{w}}

\newcommand{\x}{\bold{x}}
\newcommand{\xK}{\bold{x}_K}
\newcommand{\xs}{\bold{x}_{\s}}

\newcommand{\<}{\leqslant}

\newtheorem{Def}{Definition}[section]
\newtheorem{Lem}{Lemma}[section]
\newtheorem{Th}{Theorem}[section]
\newtheorem{Rq}{Remark}[section]

\def\righteqn#1{\noalign{\hbox to\displaywidth{\hss$\DS#1$}\penalty2000\vskip\jot}}

\title{A finite volume method on general meshes \\
for a degenerate parabolic \\
convection-reaction-diffusion equation}
\author{Ophélie Angelini \footnote{EDF R\&D, 1 avenue du Général de Gaulle 92141 Clamart, France}, Konstantin Brenner \footnote{Laboratoire de Math\' ematiques,
Universit\'e de Paris-Sud 11, F-91405 Orsay Cedex, France}, Danielle Hilhorst \footnote{CNRS and Laboratoire de Math\' ematiques,
Universit\'e de Paris-Sud 11, F-91405 Orsay Cedex, France}}

\begin{document}
\pagenumbering{roman}
\maketitle
\textbf{Abstract} We propose a finite volume method on general meshes for
the discretization of a degenerate parabolic convection-reaction-diffusion equation. Equations of
this type arise in many contexts, such as the modeling of contaminant transport in porous media.
We discretize the diffusion term, which can be anisotropic and heterogeneous, via a hybrid finite volume
scheme. We construct a partially upwind scheme for the convection term. We consider a wide range of
unstructured possibly non-matching polygonal meshes
in arbitrary space dimension. The only assumption on the mesh is that the volume elements must be star-shaped.
The scheme is fully implicit in time, it is locally conservative and robust with respect to the Péclet number.
We obtain a convergence result based upon a priori estimates and the Fr\'{e}chet--Kolmogorov compactness theorem.
\pagenumbering{arabic}

\section{Introduction}

In this paper we study a finite volume method on general meshes for
degenerate parabolic convection-reaction-diffusion equations of the form
\be\label{eq:intro}
\ddt({\beta(u)}) -\dv(\Lam\grad u)+\dv(\V u)+F(u)=q.
\ee
Equations of this type arise in particular in the modeling of contaminant transport in groundwater. The unknown function
$u$ represents the concentration of the species, which diffuses and is transported by the groundwater. An essential element in our study
is the processus of adsorbtion by a porous skeleton, which is supposed to be very fast. More particularly we suppose that the dissolved
and the absorbed parts of the species are in equilibrium; this is modeled by the function $\beta$, where $\beta'$ may be infinite in several points. The matrix $\Lam$ is a possibly anisotropic and heterogeneous diffusion-dispersion tensor, $\V$ is the velocity field, the function $F$ stands for the chemical reactions, and $q$ is the source term. We suppose that the mesh is quite general, and possibly nonmatching. Therefore, also in view of the anisotropy in the diffusion term, we can not apply the standard finite volume method \cite{Eym_Gal_Her_1}.

Finite volume schemes have often been applied to the equation (\ref{eq:intro}), see e.g. \cite{Afif_Amaz}, \cite{Baug_Walk}, \cite{Eym_Gal_Her_Mich}. The upwind discretization of the convection term permits finite volume schemes to be stable in convection dominated case, however standard finite volume schemes do not permit to handle anisotropic diffusion on general meshes. On the other hand finite element method allows a very simple discretization of full diffusion tensors, they were used a lot for the discretization of equation (\ref{eq:intro}), see e.g. \cite{Arbo_Whee_Zhan}, \cite{Daws}, \cite{Daws_Aizi}. A possible solution is to split equation (\ref{eq:intro}) into a hyperbolic part and a parabolic part, by means of an operator splitting method; one can find such an analysis in \cite{Kacu1}, \cite{Kacu2}, where the advection term was treated by the method of characteristics. The other quite intuitive idea is to take "best from both worlds" \cite{Eym_Hilh_Vohr1}, which leads to combined finite volume-finite element schemes; we refer to \cite{Eym_Hilh_Vohr1} for this approach. In order to solve this class of equations, Eymard, Hilhorst and Vohral\'{\i}k \cite{Eym_Hilh_Vohr1} discretize the diffusion term by means of piecewise linear nonconforming (Crouzeix--Raviart) finite elements over a triangularization of the space domain, or using the stiffness matrix of the hybridization of the lowest order Raviart--Thomas mixed finite element method. The other terms are discretized by means of a finite volume scheme on a dual mesh, where the dual volumes are constructed around the sides of the original triangularization. In the second paper of Eymard et al. \cite{Eym_Hilh_Vohr2} the time evolution, convection, reaction, and sources terms are discretized on a given grid, which can be nonmatching and can contain nonconvex elements, by means of a cell-centered finite volume method. In order to discretize the diffusion term, they construct a conforming simplicial mesh with vertices given by the original grid and use the finite element method. In this way, the scheme is fully consistent and the discrete solution is naturally continuous across the interfaces between the subdomains with nonmatching grids, without introducing supplementary equations and unknowns nor interpolating the discrete solutions at the interfaces.

The finite volume methods for the discretization of anisotropic diffusion on general meshes is a subject of wide interest (see for instance the results of the benchmarks organized at the FVCA 5 conference [FVCA5]). We refer to \cite{droniou_eymard_gallouet_herbin} for a detailed analysis of three recently developed families of schemes, namely the Mimetic Finite Difference scheme, the Hybrid Finite Volume scheme and the Mixed Finite Volume, which turn out to be quite similar. The most important feature of these methods is their accurate approximation of anisotropic diffusion even in highly heterogenous cases. In this paper, we apply a recent method based upon the finite volume method on general meshes developed by Eymard, Gallou\"et et Herbin~\cite{Eym_Gal_Her_1}, whereas we use a slightly modified upwind scheme for the approximation of the convection term. The time discretization is based upon a completely implicit finite difference scheme.

The organization of this paper is as follows. We describe the numerical
scheme in Section 2. We show the existence and uniqueness of the solution of the discrete scheme
and prove a priori estimates for the discrete solution
in $L^\infty(0,T;\LO)$ and in a discrete space analogous to the space
$L^2(0,T; H^1(\Om))$ in Section 3. In Section 4, we prove an estimate on differences
of time translates whereas we establish an estimate on differences of space translates
in Section 5. These estimates imply a relative compactness property of sequences of approximate
solutions by the Fr\'{e}chet--Kolmogorov theorem.  We deduce the strong convergence in $L^2$ of the
approximate solutions to the unique solution of the continuous problem in Section 6. For the proofs,
we apply methods inspired upon those of ~\cite{Eym_Gal_Her_1} and
\cite{Eym_Gal_Her_2}. In Section 7, we finally present results of numerical tests, which
confirm the validity of the numerical method.

\section{The numerical scheme}
We consider the parabolic degenerate convection-diffusion-reaction problem
$$\Pb\left\{\ba{l}
\ba{lr}
\ddt({\beta(u)}) -\dv(\Lam(\x)\grad u)+\dv(\bold{V(\x)}u)+F(u)=q(\x,t), &
(\x,t)\in Q_T,
\ea\\\\
\ba{lcr}
u(\x,t)=0 & \x\in\partial\Om, & t\in(0,T),
\ea\\\\
\ba{lcr}
u(\x,0)=u_0(\x), & \x\in\Om , &
\ea
\ea\right.$$
where $\Om \text{~is a bounded open connected polyhedral subset of~} \real^d,~d\in\intpos $, $T>0$ and $Q_T=\Om\times (0,T)$. We suppose that the following
hypotheses are satisfied:

$(\hyp_1)~~\beta\in C(\real), \beta(0)=0$ is a strictly increasing function,
which satisfies the growth condition $|\beta(a)-\beta(b)|\geqslant \underline{\beta}|a-b|, ~\underline{\beta}>0$ for all $a,b\in\real$; moreover there exist $P>0$ and $C_{\beta}$, such that $|\beta(u)|\< C_{\beta}$
for $|u|\< P$ and $\beta$ is a \Lip continuous with a constant $\overline{\beta}$ for $|u|\geqslant P$;\\
$(\hyp_2)~~\Lam$ is a measurable function from $\Om$ to $\M_d (\real)$, where  $\M_d(\real)$ denotes the set of $d\times d$ symmetric matrices, such that for a.e. $\x\in\Om$ the set of its eigenvalues is included in
$[\lambda_m,\lambda_M]$, where $\lambda_m,\lambda_M\in
L^{\infty}(\Om) \text{~are such that~} 0<\underline{\lambda}\<\lambda_m(\x)\<\lambda_M(\x)\<\overline{\lambda}$;\\
$(\hyp_3)~~\V\in H(\rm div,\Om)\bigcap L^{\infty}(\Om)$ is such that $\dv\V\geqslant 0$ a.e. in $\Om$;\\
$(\hyp_4)~~u_0\in L^{\infty}(\Om);$\\
$(\hyp_5)~~F\in C(\real)$, $F(0)=0$ and there exists $M>0$ such that $u F(u)>0$ and $F(u)$ is \Lip continuous with constant $L_{F}$ for all $u<0$ or $u>M$; moreover we suppose that $F$ does not decrease too fast i.e. there exists $\underline{F}>0$ such that $(F(u)-F(v))(u-v)\geqslant-\underline{F}(u-v)^2$ for all $u,v\in\real$;\\
$(\hyp_6)~~q\in \LQT$.\\
We now present a definition of a weak solution of Problem
$\Pb$.

\begin{Def}\label{def:weak solution}
We say that a function $u$ is a weak solution of Problem
$\Pb$ if\\
(i)~~$u\in L^2(0,T;H^1_0(\Om));$\\
(ii)~~$\beta(u)\in L^{\infty}(0,T;\LO);$\\
(iii)~~$u$~satisfies the integral equality~
$$\stl{-\ITO
\beta(u)\ph_t\dx dt-\IO
\beta(u_0)\ph(\cdot,0)\dx}
\stl{+\ITO \Lam\grad u \cdot
\grad
\ph \dx dt}$$
$$
\stl{-\ITO u\V\cdot\grad{\ph}\dx dt+\ITO
F(u)\ph\dx dt=\ITO q\ph}\dx dt
$$
\\
for all $\ph\in L^2(0,T;H^1_0(\Om))$ with
$\ph_t\in L^{\infty}(Q_T),~\ph(\cdot,T)=0$.
\end{Def}
\begin{Rq}\label{rq:unique}
In the case that the reaction function $F$ is nondecreasing, the uniqueness
of the weak solution of Problem $\Pb$ follows
from \cite{Knab_Otto}.
\end{Rq}
In order to describe the numerical scheme we introduce below some
notations related to the space and time discretization.
\begin{Def}\label{Mesh}
\textbf{(Space discretization)} Let $\Om$ be a polyhedral open
bounded connected subset of $\real^d$, with $d \in
\intpos $, and
$\partial\Om=\overline{\Om}\backslash\Om$ its boundary. A
discretization of $\Om$, denoted by $\D $, is defined
as the triplet
$\D =(\M,\E,\mathcal{P})$, where:
\par
1. $\M$ is a finite family of non empty connex open
disjoint subsets of $\Om$ (the "control volumes") such that
$\bar{\Om}=\bigcup_{K\in\M} K$. For any
$K\in\M$, let $\partial K = \overline{K}\backslash K$ be
the boundary of $K$; we define $m(K)>0$ as the measure of $K$ and
$h_K$ as the diameter of $K$.
\par
2. $\E$ is a finite family of disjoint subsets of
$\bar{\Om}$ (the "edges" of the mesh), such that, for all
$\s\in\E$, $\s$ is a non empty open subset of a
hyperplane of $\real^d$, whose $(d-1)$-dimensional measure
$m(\s)$ is strictly positive. We also assume that, for all
$K\in\M$, there exists a subset $\E_K$ of
$\E$ such that $\partial K =
\bigcup_{\s\in\E_K}\s$. For each
$\s\in\E$, we set $\M_{\s} =
\{K\in\M| \s\in\E_K\}$. We then assume that,
for all $\s\in\E$, either $\M_\s$ has
exactly one element and then $\s\in\partial\Om$ (the set of
these interfaces called boundary interfaces, is denoted by
$\E_{ext}$) or $\M_\s$  has exactly two
elements (the set of these interfaces called interior interfaces,
is denoted by $\E_{int}$). For all
$\s\in\E$, we denote by $\x_{\s}$ the
barycenter of $\s$. For all $K\in\M$ and
$\s\in\E_K$, we denote by $\n_{K,\s}$ the
outward normal unit vector.
\par
3. $\mathcal{P}$ is a family of points of $\Om$ indexed by
$\M$, denoted by $\mathcal{P} =
(\x_K)_{K\in\M}$, such that for all
$K\in\M$, $\x_K\in K$; moreover $K$ is assumed to
be $\x_K$-star-shaped, which means that for all $\x\in
K$, there holds $[\x_K, \x]\in K$. Denoting by
$d_{K,\s}$ the Euclidean distance between $\x_K$ and the
hyperplane containing $\s$, one assumes that $d_{K,\s}
> 0$. We denote by $D_{K,\s}$ the cone of vertex
$\x_K$ and basis $\s$.
\end{Def}
Next we introduce some extra notations related to the mesh. The size of
the discretization $\D $ is defined by
\be\label{def:h}
\hD=\sup_{K\in\M}{diam(K)};
\ee
moreover we define
\be\label{def:theta}
\theta_{\D}=\max(\max_{\s\in\E_{int},\{K,L\}=\mathcal{M_{\s}}}\frac{d_{K,\s}}{d_{L,\s}},
\max_{K\in\M_{\s},\s\in\E_{K}}\frac{h_{K}}{d_{K,\s}}).
\ee
Thus imposing a uniform bound on $\theta_{\D }$ forces the
meshes to be sufficiently regular. As it was done in
\cite{Eym_Gal_Her_2} we associate with the mesh the following
spaces of discrete unknowns
\be\label{X_D}
\ba{c}
X_{\D }=\{((v_K)_{K\in
\M},(v_{\s})_{\s\in \E}),
v_K\in\real,v_{\s}\in\real\},
\\ X_{\D ,0}=\{v\in X_{\D } $ such that $ (v_{\s})_
{\s\in\E_{ext}}=0\}.
\ea
\ee
Moreover, for each function $\ph=\ph(\x)$ smooth enough we define
$P_{\D}\ph\in X_D$ in following way
$$
\ba{lcl}\label{def:P_h}
(P_{\D}\ph)_{K}=\ph(\x_{K}) & $ for all $ & K\in\M,\\
(P_{\D}\ph)_{\s}=\ph(\x_{\s}) &$ for all $ &
\s\in\E.
\ea
$$
\begin{Def}\label{def:time}
\textbf{(Time discretization)}
We divide the time interval $(0,T)$ into $N$ equal time
steps of length $\dt=T/N$ such that
\be\label{hyp:k}
\dt<\underline{\beta}/\underline{F},
\ee
where $\dt$ is the uniform time step defined by $\dt=t_n-t_{n-1}$.
\end{Def}
\begin{Rq}
For the sake of simplicity, we restrict our study to the case of constant time steps.
Nevertheless all results presented below can be easily extended to the case of a non uniform time discretization.
\end{Rq}
After formally integrating the first equation of $\Pb$ on the
domain $K\times(t_{n-1},t_{n})$ for each $K\in\M$ and
$n=1,\ldots , N$, we obtain
$$
\ba{c}
\stl \int_K\beta(u(\x,t_n))-\beta(u(\x,t^{n-1}))\dx
+ \Ss \int^{t_{n}}_{t_{n-1}}\int_{\s}
(-\Lam\grad
u+\V u)\cdot\n_{K,\s}~d\gamma dt
\\
\stl
+\int^{t_{n}}_{t_{n-1}}\IK
F(u)\dx dt=\int^{t_{n}}_{t_{n-1}}\IK q\dx dt.
\ea
$$
For all $K\in\M$ and all $\s\in\E_K$ we define
$V_{K,\s}=\stl{\int_{\s}\V\cdot\n_{K,\s}}d\gamma$
and
$q^n_{K}=\stl{\frac{1}{\dt~m(K)}\int^{t_{n}}_{t_{n-1}}\IK q \dx dt}$.
We use an upwind scheme in order to approximate the
convective term, since it can possibly dominate the diffusion term; the diffusive flux
$\stl{-\int_{\s} \Lam\grad
u\cdot\n_{K,\s}d\gamma}$ is approximated by a
function of the form $F_{K,\s}(u^n)$, where $u^n=((u^n_K)_{K\in\M},
(u^n_\s)_{\s\in\E})$, and where the numerical flux $F_{K,\s}(u^n)$ is
defined by formula (\ref{def:Flux}) below. The time implicit
finite volume scheme corresponding to Problem
$\Pb$ is given by:
\par(i) The initial condition
\be\label{eq:disc:ini}
u^0_{K}=\frac{1}{m(K)}\int_K u_0(\x)\dx,
\ee
for all $K\in\M$.
\par(ii) The discrete equations
\be\label{eq:disc}
\ba{c}
\stl{m(K)(\beta(u^n_{K})-\beta(u^{n-1}_{K}))}
+\stl{\dt\Ss  {F_{K,\s}(u^n) }
+\dt\Ss  V_{K,\s}\ubar }\\+\dt~m(K)F(u^n_K)=\dt~m(K)q^n_K,
\ea
\ee
for all $K\in\M$. Unlike in the case of the standard upwind scheme, we define the value $\ubar$ as follows. For all $K\in\M$ and $\s\in\E_{K}$ we set
\be\label{upwind_int}
\ubar =\left\{ \ba{ll} u^n_K, & $ if $
V_{K,\s}\geqslant 0
\\ u^n_{\s}, & $ if $V_{K,\s}<0.
\ea\right.
\ee
We also define
\bea\label{upwind_V}
V^+_{K,\s} = \half(V_{K,\s}+|V_{K,\s}|) & \text{and} & V^-_{K,\s} = \half(V_{K,\s}-|V_{K,\s}|)
\eea
which lead to
\be\label{upwind_Vu}
V_{K,\s}\ubar =  V^+_{K,\s}u^n_K+V^-_{K,\s}u^n_{\s}
\ee
The definition of (\ref{upwind_int}) seems natural
since we also take the unknowns associated with the mesh faces.
It has an important advantage that the unknowns in the equation
(\ref{eq:disc}) are associated with a single control volume (see Remark \ref{rq:upwind});
moreover the numerical experiments presented in Section 7 show that the upwind scheme (\ref{upwind_int})
also preserves the approximate solution from unphysical oscillations
in the convection dominated case.
Finally, we remark that for each time step the number of equations is
$card(\M)$, whereas the number of discrete unknowns is
equal to $card(\M)+card(\E)$. Therefore we need to
introduce $card(\E)$ additional equations corresponding
to the interface values. For boundary faces these equations are
obtained by writing the discrete analog of the Dirichlet boundary condition
\be
\ba{lccr}
\text{(iii) }& u^n_{\s}=0 & \text{ for all }&
\s\in\E_{ext}.
\ea
\ee
For interior faces, we follow the main idea of the finite volume method by imposing the local
conservation of the discrete fluxes
\be\label{eq:conservation}
\ba{lc}
\text{(iv) }& (F_{K,\s}(u^n)+V_{K,\s}\ubar )+(F_{L,\s}(u^n)+V_{L,\s}\overline{u^n_{L,\s}})=0
\ea
\ee
for all $\s\in\E_{int}$ with $\M_{\s}=\{K,L\}$. We will define below
$F_{K,\s}$ in some more detail, but we first give an alternative variational
formulation of the discrete scheme (i)-(iv). Let $\{v^n\}_{n\in\mathbb{N}}$ be an
arbitrary sequence of elements of $X_{\D ,0}$; multiplying equation
$(\ref{eq:disc})$ by $v^n_K$ and summing on all control volumes
$K\in\M$ leads to:
$$
\ba{c}
\stl{\SK m(K)v^n_{K}\frac{\beta(u^n_{K})-\beta(u^{n-1}_{K})}{\dt}
+\SK \Ss  (v^n_{K}F_{K,\s}(u^n)+v^n_K V_{K,\s}\ubar )}
\\+\stl{\SK m(K)v^n_K F(u^n_K)=\SK m(K)v^n_K q^n_K}.
\ea
$$
Using (\ref{eq:conservation}), we obtain that
\be\label{conservativity}
\SK \Ss  v^n_{\s}(F_{K,\s}(u^n)+V_{K,\s}\ubar )=0\tfa {v^n\in X_{\D}},
\ee
which yields the following discrete weak formulation:
\par\emph{Let $u^0_{K}$ be defined by:}
\be\label{ini:var}
\ba{lr}
u^0_{K}=\stl{\frac{1}{m(K)}\int_K u_0(\x)}\dx& $
\emph{for all} $ K\in\M
\ea
\ee
\emph{For each $n\in\{1,\ldots , N\}$ find $u^n\in X_{\D ,0}$
such that for all $v^n\in X_{\D ,0}$}:
\be\label{eq:var}
\ba{c} \stl{\SK m(K)v^n_K\frac{\beta(u^n_{K})-\beta(u^{n-1}_{K})}{\dt}}+<v^n,u^n>_F+<v^n,u^n>_T\\+\stl{\SK m(K)v^n_K F(u^n_K)=\SK m(K)v^n_K q^n_K},
\ea
\ee
\emph{with}
\be\label{bil_F:var}
<v,u>_F=\SK \Ss  (v_{K}-v_{\s})F_{K,\s}(u)
\ee
\emph{and}
\be\label{bil_T:var}
<v,u>_T=\SK \Ss  (v_{K}-v_{\s})V_{K,\s}\overline{u_{K,\s}}.
\ee
Remark that the problem (i)-(iv) is equivalent to (\ref{eq:var})-(\ref{bil_T:var}). Indeed, let $\delta_{i j}$ be the Kroneker symbol, by setting $v^n_{\s}=0$ for all $\s\in\E$, and $v_K'=\delta_{K K'}$ for all $K'\in\M$ and for a given $K$ one recover (ii), and setting $v_K=0$ for all $K\in\M$ and $v_{\s'}=\delta_{\s \s'}$ for all $\s'\in\E$ yields (iv). The homogeneous Dirichlet boundary condition (iii) follows from the fact that $u^n \in X_{\D ,0}$. In order to complete the numerical scheme we still have to express the discrete
flux $F_{K,\s}$ in terms of the discrete unknowns. For this
purpose we use the SUSHI scheme proposed in \cite{Eym_Gal_Her_2}:
the idea is based upon the identification of the numerical fluxes
$F_{K,\s}$ through the mesh dependent bilinear form, using the
expression of a discrete gradient. We first define
\be\label{D_K}
\ba{lcr}
\grad_K u =
\stl{\frac{1}{m(K)}\Ss  }m(\s)(u_{\s}-u_K)\n_{K,\s}&
\forall K\in\M, & \forall u\in X_{\D }.
\ea
\ee
Remark that the geometrical relation
\be\label{eq:geom1}
\Ss m(\s)\nKs (\xs-\xK)^{T}=m(K)Id
\ee
holds for each $K$. Let $\ph(\x)$ be a function, piecewise linear
on the control volumes of the mesh. In view of (\ref{eq:geom1}) one has $\grad_K P_{\D}(\ph) = \grad\ph(\x)|_{x\in K}$. We also remark that
$$
\Ss  m(\s)\n_{K,\s}=\Ss  \int_{\s}\n_{K,\s}~d\gamma=\IK \grad
1\dx=0,
$$
which means that the coefficient of $u_K$ in $(\ref{D_K})$ is
equal to zero; thus, a reconstruction of the discrete gradient
solely based on $(\ref{D_K})$ cannot lead to a coercive
discrete bilinear form in the general case. Therefore we introduce the additional
term
\be\label{D_K,s}
\grad_{K,\s}u=\grad_{K}u+R_{K,\s}u\cdot\n_{K,\s},
\ee
where
\be\label{R_K,s}
R_{K,\s}u=\frac{\alpha_K}{d_{K,\s}}(u_{\s}-u_K-\grad_{K}u\cdot
(\bold{\x}_{\s}-\bold{\x}_K)),
\ee
for some $\alpha_K>0$, which should be chosen in a suitable way. If we choose
$\alpha_K=\sqrt{d}$ for all $K\in\M$ in the simple case that $\Lambda$ is a scalar and that the mesh
that satisfies the orthogonality property $\stl{\n_{K,\s}=\frac{\x_\s-\x_K}{d_{K,\s}}}$,
we obtain the usual two point scheme. Nevertheless, it may be
useful to optimize the choice of $\alpha_K$ as it is done in
\cite{Ang_Chav_Chen_Eym}. We then define the discrete gradient
$\grad_{\D }u$ as the piecewise constant function equal
to $\grad_{K,\s}u$ in the cone $D_{K,\s}$ with vertex
$\bold{\x}_K$ and basis $\s$
$$
\grad_{\D}u|_{D_{K,\s}}=\grad_{K,\s}u.
$$
Note that the term $R_{K,\s}$ is a
second order error term, which vanishes for piecewise linear functions. Moreover, the relation (\ref{eq:geom1}) together with (\ref{R_K,s}) implies that
\be\label{eq:E(RKs)}
\sum_{\s\in\E_K}m(D_{K,\s})R_{K,\s}(u)\n_{K,\s} = 0\tfa K\in\M \ta \tfa u\in X_{\D},
\ee
which in turn implies that
$$
\int_K \grad_{\D}u\dx = m(K) \grad_K u.
$$
The discrete gradient defined above satisfies the following strong consistency property.
\begin{Lem}\label{lemma:DiscGradLInf}
Let $\D $ be a discretization of $\Om$ in sense of
Definition \ref{Mesh}, moreover let
$\theta\geqslant\theta_{\D }$ be given. Then for all $\ph\in C^2(\overline{\Om})$, there exist a positive
constant $C$ only depending on $d$, $\theta$ and $\ph$ such
that
$$
\|\grad_{\D }P_{\D}\ph-\grad\ph\|_{(L^{\infty}(\Om))^d}\<
C \hD.
$$
\end{Lem}
The proof of this Lemma is given in \cite{Eym_Gal_Her_2}. The numerical flux is implicitly defined by the relation
\be\label{bil_Flux}
<v,u>_F=\SK \Ss  (v_{K}-v_{\s})F_{K,\s}(u)=\int_\Om
\grad_\D v\cdot\Lam(\x)\grad_\D u\dx.
\ee
It can also be defined explicitly; in order to do so, we write
the discrete gradient in the form
\be\label{D_K,s:y}
\grad_{K,\s}u=\sum_{\s'\in\E_K}(u_{\s'}-u_K)\bold{y}^{\s
\s'},
\ee
where $\bold{y}^{\s \s'}$ is defined by
$$
\bold{y}^{\s \s'}=\left\{ \ba{lc}
\stl{\frac{m(\s)}{m(K)}\n_{K,\s} +
\frac{\sqrt{d}}{d_{K,\s}}(1-\frac{m(\s)}{m(K)}\n_{K,\s}\cdot
(\x_{\s}-\x_K))\n_{K,\s}} & $ if $ \s=\s',\\
\stl{\frac{m(\s')}{m(K)}\n_{K,\s'} -
\frac{\sqrt{d}}{d_{K,\s}}\frac{m(\s')}{m(K)}\n_{K,\s'}\cdot
(\x_{\s}-\x_K)\n_{K,\s}} & $ otherwise$.
\ea\right.
$$
We obtain that for all $u,v\in X_{\D }$
$$
\IO \grad_{\D }u(\x)\cdot
\Lam(\x)\grad_{\D }v(\x)=\SK \Ss \sum_{\s'\in\E_K}
A^{\s \s'}_{K}(u_{\s}-u_K)(v_{\s'}-v_K),
$$
with
$$
A^{\s \s'}_{K}=\stl{\sum_{\s''\in\E_K}\bold{y}^{\s''
\s}\cdot \Lam_{K,\s''}\bold{y}^{\s'' \s'}}
\ta
\Lam_{K,\s''}=\stl{\int_{D_{K,\s''}}\Lam(\x)}\dx.
$$
The local matrices $(A^{\s \s'}_{K})_{\s
\s'\in\E_K}$ are symmetric, and the numerical flux is
then defined by
\be\label{def:Flux}
F_{K,\s}(u)=\sum_{\s'\in\E_K}A^{\s
\s'}_{K}(u_{K}-u_{\s'}).
\ee
Next we prove some useful properties of the
mesh depending bilinear forms introduced previously. In
particular we prove that $<\cdot,\cdot>_F$ is continuous and
coercive, and that $<\cdot,\cdot>_T$ is continuous and
nonnegative. The space $X_{\D }$ defined in (\ref{X_D})
is equipped with the following semi-norm.
\begin{Def}
Let $\D =(\M,\E,\mathcal{P})$ be a
discretization of $\Om$ in the sense of Definition \ref{Mesh};
then for all $v\in X_{\D }$ we define
\be
|v|^2_{X}=\SK \Ss  \frac{m(\s)}{d_{K,\s}}(v_{\s}-v_K)^2,
\ee
\end{Def}
which is a norm on the space $X_{\D ,0}$. Let us also define a discrete analog of $\|\cdot\|_{1,p}$ norm.
\begin{Def}\textbf{(The discrete space $H_\M(\Om)$)}\label{def:W1p} Let $1\< p<\infty$ and let $\D =(\M,\E,\mathcal{P})$ be a discretization of $\Om$ in the sense of
Definition \ref{Mesh}. Let $H_\M(\Om)\subset
\LO$ be the set of piecewise constant functions on the
control volumes of the mesh $\M$ for each $v\in H_\M(\Om)$ we define $v_K=v(\x)|_{\x\in K}$.
\\
For all $v\in
H_\M(\Om)$ and for all $\s\in\E_{int}$
with $\M_{\s}=\{K,L\}$ we define
$D_{\s}v=|v_K-v_L|$ and
$d_{\s}=d_{K,\s}+d_{L,\s}$, and for all
$\s\in\E_{ext}$ with $\M_{\s}=\{K\}$, we
set $D_{\s}v=|v_K|$ and $d_{\s}=d_{K,\s}$. We then
define the following family of norms
\be\label{def:W1p_eq}
\|v\|^p_{1,p,\M}=\SK \Ss  m(\s)d_{K,\s}
(\frac{D_\s v}{d_{\s}})^p;
\ee
\end{Def}
so that in particular
$$
\|v\|^2_{1,2,\M}=\SK \Ss  m(\s)d_{K,\s}(\frac{D_{\s}v}{d_{\s}})^2.
$$
Next we recall two results from \cite{Eym_Gal_Her_2} which
we will use below. The following lemma shows the equivalence between the semi-norm in
$X_{\D }$ and the $L^2$-norm of the discrete gradient.
\begin{Lem}\label{lemma:equivalence}
Let $\D $ be a discretization of $\Om$ in the sense of
Definition \ref{Mesh}, and let
$\theta\geqslant\theta_{\D }$ be given. Then there exists
$C_1>0$ and $C_2>0$ only depending on $\theta$ and $d$ such that
$$
\ba{lr}
C_1 |v|_{X}\<
\|\grad_\D v\|_{\LO}\< C_2 |v|_{X} & $for
all $ v\in X_{\D }.
\ea
$$
\end{Lem}
\begin{Lem}\label{lemma:w12bound}
Let $\D $ be a discretization of $\Om$ in the sense of
Definition \ref{Mesh}, then there holds
$$
\ba{lr}
\|v\|_{1,2,\M}\< |v|_{X} & $for all $ v\in
X_{\D ,0}.
\ea
$$
\end{Lem}
Next we show that the bilinear forms defined in
(\ref{bil_F:var}) and (\ref{bil_T:var}) satisfy continuity and coercivity
properties.
\begin{Lem}\label{lemma:bils}
Let $\D $ be a discretization of $\Om$ in the sense of
Definition \ref{Mesh}, and let
$\theta\geqslant\theta_{\D }$ be given, then:
\par(i) There exist positive constants $C_1$ and $\alpha$ which
do not depend on $h$ such that
$$
|<u,v>_F|\< C_1 |u|_X |v|_X
$$
and
$$
<u,u>_F\geqslant \alpha |u|^2_X
$$
for all  $u,v\in X_{\D }$.
\par (ii) There exist a
positive constant $C_2$ which does not depend on $h$
that
$$
|<u,v>_T|\< C_2 |u|_X |v|_X
$$
and
$$
<u,u>_T\geqslant 0
$$
for all  $u,v\in X_{\D ,0}$.

\end{Lem}
\emph{Proof.} (i) Using the definition of the numerical flux
(\ref{bil_Flux}) and in view of $(\hyp_2)$ and Lemma \ref{lemma:equivalence}
$$
|<u,v>_F|=|\int_\Om\grad_\D u\cdot\Lam(x)\grad_\D v|\dx\<
\overline{\lambda}\|\grad_\D u\|_{\LO}\|\grad_\D v\|_{\LO}
\< C_1|u|_X |v|_X;
$$
on the other hand we have that
$$
<u,u>_F=\int_\Om
\grad_\D u\cdot\Lam(x)\grad_\D u \dx\geqslant
\underline{\lambda}\dx\|\grad_\D u\|^2_{\LO}
\geqslant C_2|u|^2_X.
$$
\par (ii) By the definition (\ref{bil_T:var}) we have that
$$
<u,v>_T=\SK \Ss  (v_K-v_{\s})
V_{K,\s} \overline{u_{K,\s}}.
$$
Using the definition (\ref{upwind_int}) one can write:
$$
<u,v>_T=\stl\SK \sum_{\s\in{\E_K},V_{K,\s}\geqslant 0} V_{K,\s}(v_K-v_{\s})u_K +
\SK \sum_{\s\in{\E_K},V_{K,\s}\< 0}V_{K,\s}(v_K-v_{\s})u_{\s},
$$
which implies
\be\label{loc:bil_T}
<u,v>_T=\stl\SK \Ss  V_{K,\s}(v_K-v_{\s})u_K-
\SK \sum_{\s\in{\E_K},V_{K,\s}\< 0}V_{K,\s}(v_K-v_{\s})(u_K-u_{\s}).
\ee
Using the Cauchy-Schwarz inequality and the bound $d_{K,\s}\< \hD$ we have that
$$
\bl
|<u,v>_T|&
\<\stl{\sqrt{d}\cdot\Vmax(\SK \Ss
m(\s)\frac{(v_K-v_{\s})^2}{d_{K,\s}})^{\half}(\SK \Ss  \frac{m(\s)d_{K,\s}}{d}u^2_K)^{\half}}\\
&+\hD\cdot\stl{\Vmax (\SK \Ss
m(\s)\frac{(v_K-v_{\s})^2}{d_{K,\s}})^{\half}(\SK \Ss
m(\s)\frac{(u_K-u_{\s})^2}{d_{K,\s}})^{\half}}.
\el
$$
and
$$
|<u,v>_T|\< \Vmax (\sqrt{d}\cdot |v|_{X}\|u\|_{\LO}+diam(\Om)\cdot|v|_{X}|u|_{X})
$$
since $\stl{\sum_{\s\in{\E_{K}}}d_{K
\s}} m(\s)=m(K)d$.
In view of Lemma \ref{lemma:w12bound} and the discrete Poincar\'{e}
inequality implied by Lemma \ref{lemma:LqW1p} below we conclude that
$$
|<u,v>_T|\< C_2|u|_X |v|_X.
$$
In order to prove the positivity, we write $<u,u>_T$ in the form (\ref{loc:bil_T})
$$
<u,u>_T=
\stl\SK \Ss  V_{K,\s}(u_K-u_{\s})u_K-
\SK \sum_{\s\in{\E_K},V_{K,\s}\< 0}V_{K,\s}(u_K-u_{\s})^2;
$$
using the algebraic inequality $-2ab\geqslant -a^2-b^2$ and the discrete boundary condition we obtain
$$
<u,u>_T\geqslant \stl\half\SK \Ss V_{K,\s} (u^2_K-u^2_{\s})=
\stl{\half\SK u^2_K\Ss  V_{K,\s}}.
$$
By the assumption ($\hyp_3$) one has that
$\stl{\sum_{\s\in{\E_{K}}}V_{K,\s}\geqslant 0}$ and we finally conclude that
$$
<u,u>_T\geqslant 0.
$$
\\
Next we recall a technical lemma presented in \cite{Eym_Hilh_Vohr1}, Lemma 8.2,
which will be useful for the a priori estimates of the next section
\begin{Lem}\label{lemma:B}
Let $B(s)$, $s\in \real$ be defined by
$$
\stl{B(s)=\beta(s)s-\int^s_0\beta(\tau)d\tau},
$$
with $\beta$ satisfying hypothesis $(\hyp_1)$. Then $\stl{B(s)\geqslant \half s^2 \underline{\beta}}$.
\end{Lem}

\section{A priori estimates.}
We define below an approximate solution of Problem (1)-(3).
\begin{Def}\label{def:solution}\textbf{(Approximate solution)}\\
Let the sequence of $\{u^n\}\in X^N_{\D ,0}$, $n\in\{1,\ldots, N\}$, be
a solution of the discrete problem (\ref{ini:var})-(\ref{bil_T:var}), with $\dt=T/N>0$.
We say that the piecewise constant function $\uhk:\Om\times [0,T]\to\real$ is an
approximate solution of Problem
$\Pb$ if
$$
\ba{ll}
\uhk(\x,0)=u^0_K & $ for all $\x\in K,\\
\uhk(\x,t)=u^n_K & $ for all $(\x,t)\in K\times(t_{n-1},t_n];\\
\ea
$$
we also define its approximate gradient by
$$
\ba{rcl}
\grad_{\D,\dt} \uhk(\x,t)=\grad_{\D }u^n(\x)& $ for all $(\x,t)\in K\times(t_{n-1},t_n].\\
\ea
$$
\end{Def}
\begin{Lem}\textbf{(A priori estimate)}
Let $\uhk$ be an approximate solution of Problem (1)-(3), then it is such that
\be\label{ineq:energy}
\stl{\frac{1}{4}\underline{\beta}\|\uhk\|^2_{L^{\infty}(0,T;\LO)} \< C_1 ~\ta~ \stl{\frac{1}{2}\underline{\lambda}\|\grad_{\D,\dt} \uhk\|^2_{\LQT}}\< C_1},
\ee
where
$$
\stl{C_1=\frac{1}{\underline{\beta}}\|\beta(u_0)\|^2_{L^2(\Om)}+m(\Om)TM |\min_{0\< u\< M}F(u)|+\frac{T}{\underline{\beta}}\|q\|^2_{\LQT}};
$$
moreover there exists $C_2>0$, such that
\be\label{ineq:b(u)}
\|\beta(\uhk)\|_{L^{\infty}(0,T;\LO)} \< C_2.
\ee
\end{Lem}
\emph{Proof.} Let $m \in [1,N]$ be an arbitrary integer. Summing on $n \in \{1,\ldots , m\}$ the equation (\ref{eq:var})
with $v^n=u^n$ for each $n$ we obtain
$$
\ba{c}\stl{
\SK m(K)\sum^m_{n=1}u^n_{K}(\beta(u^n_{K})-\beta(u^{n-1}_{K}))+\sum^m_{n=1}\dt(<u^{n},u^{n}>_F+<u^{n},u^{n}>_T)}\\
\stl{+\sum^m_{n=1}\dt\SK m(K)u^n_K F(u^n_K)=\sum^m_{n=1}\SK \dt~m(K)u^n_K q^n_K.
}\ea
$$
Next, we consider the function $B$ from Lemma \ref{lemma:B} defined by
$$
B(u)=\beta(u)u-\int^{u}_0 \beta(\tau)d\tau.
$$
One can see that the following relation holds
$$
B(u^n_K)-B(u^{n-1}_K)=u^n_K(\beta(u^n_K)-\beta(u^{n-1}_K))-\int^{u^n_K}_{u^{n-1}_K} (\beta(\tau)-\beta(u^{n-1}_K))d\tau
$$
and since $\beta$ is nondecreasing we have that
$$
\int^{u^n_K}_{u^{n-1}_K} (\beta(\tau)-\beta(u^{n-1}_K))d\tau\geqslant 0,
$$
which implies
$$
\bl\stl
\SK m(K)(B(u^m_{K})-B(u^0_{K}))
&=\stl\SK m(K)\sum^m_{n=1}(B(u^n_{K})-B(u^{n-1}_{K}))\\
&\stl\<\SK m(K)\sum^m_{n=1}u^n_{K}(\beta(u^n_{K})-\beta(u^{n-1}_{K})).
\el
$$
In view of Lemma \ref{lemma:B} we have that
$$
\half\underline{\beta}u^2\< B(u)\< u\beta(u)\< \frac{(\beta(u))^2}{\underline{\beta}},
$$
which yields
$$
\stl{
\half\underline{\beta}\|\uhk(\cdot,t_m)\|^2_{\LO}-\frac{1}{\underline{\beta}}\|\beta(u_0)\|^2_{L^2(\Om)}}\<\SK m(K)\sum^m_{n=1}u^n_{K}(\beta(u^n_{K})-\beta(u^{n-1}_{K})).
$$
We remark that in view of the hypothesis $(\hyp_5)$ one has
$$u F(u)\geqslant \stl{M\min_{0\< u\< M}F(u)},$$
since $\stl{\min_{0\< u\< M}F(u)}\< 0$. The last statement of Lemma \ref{lemma:bils} implies $<u^n,u^n>_T\geqslant 0 $. By the equation (\ref{bil_Flux}) and $(\hyp_2)$ we finally conclude that
\be\label{loc:a_priori}
\ba{c}
\stl
\frac{1}{2}\underline{\beta}\|\uhk(\cdot,t_m)\|^2_{\LO}+\underline{\lambda}\|\grad_{\D,\dt}
\uhk\|^2_{L^2(\Om\times(0,k m))}\\
\<\stl C+\sum^m_{n=1}\SK \dt m(K)u^n_K q^n_K,
\ea
\ee
where
$$
\stl{C=\frac{1}{\underline{\beta}}\|\beta(u_0)\|^2_{L^2(\Om)}+m(\Om)TM |\min_{0\< u\< M}F(u)|}.
$$
Applying Cauchy-Schwarz and Young's inequality to the last term in (\ref{loc:a_priori}) leads to
$$
\ba{c}
\stl{\frac{1}{2}\underline{\beta}\|\uhk(\cdot,t_m)\|^2_{\LO}+\underline{\lambda}\|\grad_{\D,\dt}
\uhk\|^2_{L^2(\Om\times(0,k m))}\<  C+\|\uhk\|_{\LQT}\|q\|_{\LQT}}\\
\< \stl{ C+\frac{\varepsilon}{2}\|\uhk\|^2_{\LQT}+\frac{1}{2\varepsilon}\|q\|^2_{\LQT}}.
\ea
$$
We then obtain
$$
\ba{c}
\stl{\half\underline{\beta}\|\uhk\|^2_{L^{\infty}(0,T;\LO)} \<  C+\frac{\varepsilon}{2}T\|\uhk\|^2_{L^{\infty}(0,T;\LO)}+\frac{1}{2\varepsilon}\|q\|^2_{\LQT}}
\ea
$$
and
$$
\ba{c}
\stl{\underline{\lambda}\|\grad_{\D,\dt} \uhk\|^2_{\LQT}} \<  \stl{C+\frac{\varepsilon}{2}T\|\uhk\|^2_{L^{\infty}(0,T;\LO)}+\frac{1}{2\varepsilon}\|q\|^2_{\LQT}}.
\ea
$$
We now choose $\varepsilon=\underline{\beta}/(2T)$, which gives
$$
\stl{\frac{1}{4}\underline{\beta}\|\uhk\|^2_{L^{\infty}(0,T;\LO)}} \<  C ~\ta~ \stl{\frac{\underline{\lambda}}{2}\|\grad_{\D,\dt} \uhk\|^2_{\LQT}}\< C,
$$
where
$$
\stl{C=\frac{1}{\underline{\beta}}\|\beta(u_0)\|^2_{L^2(\Om)}+m(\Om)TM |\min_{0\< u\< M}F(u)|+\frac{T}{\underline{\beta}}\|q\|^2_{\LQT}}.
$$
In order to prove the estimate on $\|\beta(\uhk)\|_{L^{\infty}(0,T;\LO)}$
we split $\beta$ into a bounded and a \Lip continuous part by setting $\beta = \beta_1+\beta_2$, where
\begin{center}
\be\label{def:splitting1}
\beta_1(s)=
\left\{\ba{cc}
\beta(s)& |s|\< P\\
0&|s|> P,
\ea\right.
~~~\beta_2(s)=
\left\{\ba{cc}
0&|s|\< P\\
\beta(s)& |s|> P,
\ea\right.
\ee
\end{center}
and
\be\label{def:splitting2}
y(s)=\left\{
\ba{cc}
\stl{\frac{\beta(P)-\beta(-P)}{2 P}s+\frac{\beta(P)+\beta(-P)}{2}}&|s|\< P\\
0&|s|>P.
\ea
\right.
\ee
We finally define
\be\label{def:splitting3}
\widetilde{\beta}_1=\beta_1-y ~\ta~ \widetilde{\beta}_2=\beta_2+y.
\ee
we then remark that $\beta=\widetilde{\beta}_1+\widetilde{\beta}_2$; we remark that $\widetilde{\beta}_1$ and $\widetilde{\beta}_2$
are continuous and that $\widetilde{\beta}_1$ is bounded by $2C_{\beta}$, while $\widetilde{\beta}_2$ is  \Lip continuous with \Lip constant $L_{\widetilde{\beta}}=\max(\overline{\beta},(\beta(P)-\beta(-P))/2 P)$. Which implies the $\LinfO$ estimate
$$
\|\beta(\uhk)\|_{\LinfO }\<\|\widetilde{\beta}_1(\uhk)\|_{\LinfO }
$$
$$
+\|\widetilde{\beta}_2(\uhk)-\widetilde{\beta}_2(0)\|_{\LinfO }+\|\widetilde{\beta}_2(0)\|_{\LinfO },
$$
so that
$$
\|\beta(\uhk)\|_{\LinfO }\<2m(\Om)^{\half}C_{\beta}+L_{\widetilde{\beta}}\|\uhk\|_{\LinfO }+m(\Om)^{\half}|\widetilde{\beta}_2(0)|.
$$
\begin{Rq}\label{rq:eps_pb}(\textbf{Extended discrete problem})
Let $s>0$ and $w\in H_{\M}(\Om)$ (sf. Definition \ref{def:W1p}), we consider the following extended one step problem.
Find $u\in X_{\D ,0}$ such that for all $v\in X_{\D ,0}$:
\be\label{loc:extended}
\ba{c} \stl{s\SK m(K)v_K\frac{\beta(u_{K})-\beta(w_{K})}{\dt}}+<v,u>_F+s< v,u>_T\\+\stl{s\SK m(K)v_K F(u_K)=s\SK m(K)v_K q_K}.
\ea
\ee
It can be shown that the solution of the extended problem (\ref{loc:extended}) is bounded in the norm $|\cdot|_{X}$. More precisely, it is such that
\be\label{ineq:energy_extended}
\stl{\frac{1}{2}\underline{\lambda}\|\grad_{\D,\dt} \uhk\|^2_{\LQT}}\< s(\frac{1}{\underline{\beta}}\|\beta(w)\|^2_{L^2(\Om)}+m(\Om)TM |\min_{0\< u\< M}F(u)|+\frac{T}{\underline{\beta}}\|q\|^2_{\LQT}).
\ee
\end{Rq}
\begin{Th}\textbf{(Existence of a discrete solution)}
The problem (\ref{ini:var})-(\ref{bil_T:var}) has at least one solution.
\end{Th}
\emph{Proof.} Let $(e_i)_{1\< i\< card(X_{\D ,0})}$ be a family elements of $X_{\D ,0}$, which components are defined by $(e_i)_j=\delta_{i j}$, where $\delta_{i j}$ is the Kronecker symbol. The system of nonlinear equations (\ref{eq:disc:ini})-(\ref{eq:conservation}) may be written in the form
\be\label{loc:ex}
E(\beta(u^n)-\beta(u^{n-1}))+\mat{A}u^n+\mat{C}u^n+\dt E F(u^n)=\mat{Q}^n,
\ee
where
\par (i) $u^{n}$, $u^{n-1}\in X_{\D ,0}$;
\par (ii) $E$ is the diagonal matrix of the size $card(\M)+card(\E_{int})$ with elements
$$
(E)_{K,K}=m(K)~~\text{and}~~(E)_{\s,\s}=0
$$
for all $K\in\M,\s\in\E_{int}$;
\par (iii) $\beta$ and $F$ are continuous mappings from $X_{\D ,0}$ to itself naturally defined by
$$
(\beta(u))_i =\beta(u_i)~~\text{and}~~(F(u))_i = F(u_i);
$$
\par (iv) $\mat{A}$ and $\mat{C}$ are the diffusion matrix and the convection matrix respectively, with components
$$
\mat{A}_{i j}=\dt<e_i,e_j>_F\ta\mat{C}_{i j}=\dt<e_i,e_j>_T,
$$
(v) $\mat{Q}^n\in X_{D,0}$ is the source term, given by
$$
\mat{Q}^n_K=\dt m(K)q^n_K~\text{and}~\mat{Q}^n_{\s}=0
$$
for all $K\in\M,\s\in\E_{int}$;\\\\
Due to the coercivity of the bilinear form corresponding to the diffusion the matrix $\mat{A}$ is invertible; hence $(\ref{loc:ex})$ is equivalent to
$$
u^n+\mat{A}^{-1}(E(\beta(u^n)-\beta(u^{n-1}))+\mat{C}u^n+\dt E F(u^n)-\mat{Q}^n)=0.
$$

As it has been done in (\ref{loc:extended}), we introduce the extended formulation
\be\label{loc:ex1}
u^n+s\mat{A}^{-1}(E(\beta(u^n)-\beta(u^{n-1}))+\mat{C}u^n+\dt E F(u^n)-\mat{Q}^n)=0,
\ee
with $s\in[0,1]$. Moreover for a given $u^{n-1}$ we define a continuous mapping $H_n:[0,1]\times X_{\D ,0}\to X_{\D ,0}$ by
$$
H_n(s,u)=s \mat{A}^{-1}(E(\beta(u)-\beta(u^{n-1}))+\mat{C}u+\dt E F(u)-\mat{Q}^n).
$$
Then the equation (\ref{loc:ex1}) can be written in the form $u^n+H_n(s,u^n)=0$. In view of Remark \ref{rq:eps_pb}, the estimate ($\ref{ineq:b(u)}$) and the Lemma \ref{lemma:equivalence} we have that
$$
\dt|u|^2_{X}\< C,
$$
with some positive constant $C$, which does not depend on $s$. Setting $R=\sqrt{C/\dt+1}$ we deduce that
$$
\ba{c}
|u|_{X}<R \tfa (s,u)\in [0,1]\times X_{\D ,0} \text{~such that~} u+H_n(s,u)=0.
\ea
$$
Therefore the equation $u+H_n(s,u)=0$ has no solutions on the boundary of the ball $B_{R}$ of radius $R$ for $s\in[0,1]$. Next, we denote by $d(Id+H_n(s,\cdot),B_R,0)$ the topological degree of the application $Id+H_n(s,\cdot)$ with respect to the ball $B_R$ and \rhs $0$. In view of the homotopy invariance of the topological degree and thanks to the fact that $H_n(0,u)=0$ for all $u\in X_{\D ,0}$ we have that
$$
d(Id+H_n(s,\cdot),B_R,0)=d(Id+H_n(0,\cdot),B_R,0)=1\tfa s\in[0,1],
$$
where we have applied [\cite{Deimling}, Theorem 3.1 (d1) and (d3)]. Thus, by [\cite{Deimling}, Theorem 3.1 (d4)], there exists $u^n$ such that $u^n+H_n(1,u^n)=0$, so that $u^n$ is a solution of (\ref{loc:ex}).
\begin{Th}\textbf{(Uniqueness of the discrete solution)}
The solution of the problem (\ref{ini:var})-(\ref{bil_T:var}) is unique.
\end{Th}
\emph{Proof.}
We give a proof by contradiction. Let $\uhk$ and $\widetilde{u}_{h,k}$ be two different solutions of (\ref{ini:var})-(\ref{bil_T:var}), such that $u^m=\widetilde{u}^m$ for all $m=1,\ldots , n-1$, but $u^n\neq\widetilde{u}^n$. We define $r^n=u^n-\widetilde{u}^n$. In view of (\ref{eq:var}) with $v=r^n$ we have that
$$
\stl{\SK m(K)r^n_K\frac{\beta(u^n_{K})-\beta(\widetilde{u}^n_{K})}{\dt}}+<r^n,r^n>_F+<r^n,r^n>_T
$$
$$
+\stl{\SK m(K)r^n_K (F(u^n_K)-F(\widetilde{u}^n_K))=0}.
$$
We apply Lemma \ref{lemma:bils} as well as the assumptions ($\hyp_1$) and ($\hyp_5$) in order to estimate each term in the above equation. We obtain that
$$
\stl{({\underline{\beta}}/{\dt}-\underline{F})\SK m(K)(r^n_K)^2}+\alpha |r^n|^2_{X}\< 0,
$$
where $\alpha$ is the coercivity constant. Finally, in view of the assumption (\ref{hyp:k}) on the time step we deduce that
$$
|r^n|_{X}=0.
$$

\section{Estimate on time translates}
To begin with we give two technical lemmas which will be useful for proving the estimate
on time translates

\begin{Lem}\label{lemma:EGH1} Let $T>0$, $\tau\in(0,T)$, $N\in\intpos$, $\dt=T/N$ be given and
$(a^n)_{n\in\intpos }$ be a family of non negative real values. Let $\lceil s \rceil$ denotes the smallest integer larger or equal to $s$.
Then
$$
\int^{T-\tau}_0\sum_{\lc t/\dt\rc+1\<n\<\lc(t+\tau)/\dt\rc} a^n dt\< \tau
\Sn a^n.
$$
\end{Lem}
\emph{Proof.}
One has that
$$
  \int^{T-\tau}_0\sum_{\lc t/\dt\rc+1\<n\<\lc(t+\tau)/\dt\rc} a^n dt
\<\int^{T-\tau}_0\sum_{t/\dt+1\<n<(t+\tau)/\dt+1} a^n dt
 =\int^{T-\tau}_0\sum_{t\<m \dt<t+\tau} a^{m+1} dt
$$
We remark that if ${\lc t/\dt\rc+1>\lc(t+\tau)/\dt\rc}$, then the above inequality seal holds, with the left hand side term equal to zero.
We define a characteristic function $\chi(n,t_1,t_2)$ by
$$\chi(n,t_1,t_2)=\left\{\ba{cc}
1 & \text{~if~} t_1\< n \dt<t_2,\\
0 & \text{otherwise.}
\ea\right.
$$
Then we obtain that
$$
\int^{T-\tau}_0\sum_{\lc t/\dt\rc+1\<n\<\lc(t+\tau)/\dt\rc} a^n dt
\< \sum^{N-1}_{m=1}a^{m+1}\int^{T-\tau}_0\chi(n,t,t+\tau)dt\<\tau
\sum^{N}_{m=1}a^m.
$$

\begin{Lem}\label{lemma:EGH2}
Let $T>0$, $\tau\in(0,T)$, $N\in\intpos$, $\dt=T/N$, $\zeta\in[0,\tau]$ be given and
$(a^n)_{n\in\intpos }$ be a family of nonnegative real values. Let $\lceil s \rceil$ denotes the smallest integer larger or equal to $s$. Then
$$
\int^{T-\tau}_0\sum_{\lc t/\dt\rc+1\<n\<\lc(t+\tau)/\dt\rc} a^{\lc (t+\zeta)/\dt\rc}dt\< \tau
\Sn a^n.
$$
\end{Lem}
\emph{Proof.}
As in the proof of the previous Lemma we have that
$$
\int^{T-\tau}_0\sum_{\lc t/\dt\rc+1\<n\<\lc(t+\tau)/\dt\rc} a^{\lc (t+\zeta)/\dt\rc}dt\<\Sn \int^{T-\tau}_0 a^{\lc (t+\zeta)/\dt\rc}\chi(n,t,t+\tau)dt
$$
A simple change of variable implies
$$
\bl
\stl \Sn \int^{T-\tau}_0 a^{\lc (t+\zeta)/\dt\rc}\chi(n,t,t+\tau)dt
&\stl \<\Sn  \int^{T}_0 a^{\lc s/\dt\rc}\chi(n,s-\zeta,s-\zeta+\tau)ds\\
&\stl =\sum^{N}_{m=1}a^m \Sn  \int^{m \dt}_{(m-1)\dt}\chi(n,s-\zeta,s-\zeta+\tau)ds.
\el
$$
In order to conclude the proof we remark that $\chi(n,t_1,t_2)=\chi(n+m,t_1+m \dt,t_2+m \dt)$ for all $n,m\in\mathbb{Z},~t_1,t_2\in\real$, which in turn implies that
$$
\bl
\stl \Sn  \int^{m \dt}_{(m-1)\dt}\chi(n,s-\zeta,s-\zeta+\tau)ds
&\stl = \Sn  \int^{-n \dt}_{-(n+1)\dt}\chi(-m,s-\zeta,s-\zeta+\tau)ds\\
&\stl \<\int_{\real}\chi(-m,s-\zeta,s-\zeta+\tau)ds= \tau.
\el
$$
\begin{Th}
Let $\D $ be a discretization of $\Om$ in the sense of
Definition \ref{Mesh} and let $\{\uhk\}$ be a solution of the discrete
problem in sense of Definition \ref{def:solution}. Let also
$\theta\geqslant\theta_{\D }$ be given. Then there exists a
positive constant $C$ only depending on $\theta$ such that
\be\label{ineq:TTrans}
\int^{T-\tau}_{0}\IO(\uhk(x,t+\tau)-\uhk(x,t))^2 dx
dt \< C\tau,
\ee
for all $\tau\in(0,T)$.
\end{Th}
\emph{Proof.}
To begin with we use the hypothesis  ($\hyp_1$) to obtain
\bean
\lqn{\underline{\beta}\int^{T-\tau}_{0}\IO (\uhk(x,t+\tau)-\uhk(x,t))^2\dx dt} & &\\
& &=\underline{\beta}\int^{T-\tau}_{0}\SK m(K)(u^{\lc (t+\tau)/\dt \rc}_K-u^{\lc t/\dt\rc}_K)^2dt\\
& &\<\int^{T-\tau}_{0}\SK m(K)(u^{\lc (t+\tau)/\dt \rc}_K-u^{\lc t/\dt\rc}_K)(\beta(u^{\lc (t+\tau)/\dt \rc}_K)-\beta(u^{\lc t/\dt\rc}_K))dt\\
& &=\int^{T-\tau}_{0}\SK m(K)(u^{{\lc (t+\tau)/\dt \rc}}_K-u^{{\lc t/\dt\rc}}_K)\sum_{\lc t/\dt \rc+1\<n\<\lc (t+\tau)/\dt \rc}m(K)(\beta(u^{n}_K)-\beta(u^{n-1}_K))dt
\eean
For a given $k$ and for all real $t$ and $\tau$ we define the following set
$$
n(t,\tau)=\{n\in\mathbb{N} ,~{\lc t/\dt \rc+1\<n\<\lc (t+\tau)/\dt \rc}\},
$$
which can be empty. Then, the discrete equation (\ref{eq:disc}) implies
\bean
\lqn{\underline{\beta}\int^{T-\tau}_{0}\IO (\uhk(x,t+\tau)-\uhk(x,t))^2\dx dt\<\int^{T-\tau}_{0}\SK (u^{\lc (t+\tau)/\dt\rc}_K-u^{\lc t/\dt\rc}_K)}\\
& & \cdot\sum_{n\in n(t,\tau)}\dt(m(K)q^{n}_K-\Ss (F_{K,\s}(u^{n})+V_{K,\s}\ubar )-m(K)F(u^{n}_K))dt.
\eean
Let us define the expressions $A_{D,C}$, $A_R$ and
$A_S$ by

\bean
\lqn{\stl
A_{D,C}=\int^{T-\tau}_{0}\sum_{n\in n(t,\tau)}\dt \SK (u^{\lc (t+\tau)/\dt\rc}_K-u^{\lc t/\dt\rc}_K)\Ss (F_{K,\s}(u^{n})+V_{K,\s}\ubar )dt,}\\
& \lqn{
\stl A_R=\int^{T-\tau}_{0}\sum_{n\in n(t,\tau)}\dt \SK m(K)(u^{\lc (t+\tau)/\dt\rc}_K-u^{\lc t/\dt\rc}_K)F(u^{n}_K)dt,}
\\
& &\stl A_S=\int^{T-\tau}_{0}\sum_{n\in n(t,\tau)}\dt\SK m(K)(u^{\lc (t+\tau)/\dt\rc}_K-u^{\lc t/\dt\rc}_K) q^{n}_K dt,
\eean
which we will estimate below. In view of (\ref{conservativity}), (\ref{bil_F:var}) and (\ref{bil_T:var})
we obtain
$$
A_{D,C}=
\int^{T-\tau}_{0}\sum_{n\in n(t,\tau)}\dt(
<u^{\lc (t+\tau)/\dt\rc}-u^{\lc t/\dt\rc},u^{n}>_{F}+<u^{\lc (t+\tau)/\dt\rc}-u^{\lc t/\dt\rc},u^{n}>_{T})dt
$$
In view of Lemma $\ref{lemma:bils}$ we have that $\stl{|<u,v>_{F}|+|<u,v>_{T}|\< C|u|_X |v|_X}$ for all $u,v\in X_{\D ,0}$ and since $2a b\< a^2+b^2$ one has
$$
|A_{D,C}|\<C\int^{T-\tau}_{0}\sum_{n\in n(t,\tau)}\dt(|u^{\lc (t+\tau)/\dt\rc}|_X+|u^{\lc t/\dt\rc}|_X) |u^{n}|_X dt
$$
$$
\< C(\int^{T-\tau}_{0}\sum_{n\in n(t,\tau)}\dt|u^{\lc t/\dt\rc}|^2_X+\int^{T-\tau}_{0}\sum_{n\in n(t,\tau)}\dt|u^{\lc (t+\tau)/\dt\rc}|^2_X+
\int^{T-\tau}_{0}\sum_{n\in n(t,\tau)}\dt|u^{n}|^2_X dt).
$$
It follows from the estimate (\ref{ineq:energy}) and the Lemmas \ref{lemma:EGH1} and \ref{lemma:EGH2}
that
$$
|A_{D,C}|\<\tau C\sum_{n=1}^{N}\dt|u^{n}|^2_X\< C\tau.
$$
Next, we estimate the term $A_R$; we remark that for all $u,v\in\real$ it holds
$$
v F(u)\<L_F|v||u|\<\half v^2 + \half L^2_F u^2\text{~~if~} u< 0,
$$
$$
\stl{v F(u)\<|v|\max_{0\<u\<M}|F(u)|\< \half v^2 + \half \max_{0\<u\<M}F^2(u)}\text{~~if~} 0\<u\<M,
$$
$$
v F(u)\<|v|(|F(u)-F(M)|+|F(M)|)\<|v|(L_F|u-M|+F(M))
$$
$$
\<v^2+\half L^2_F|u|^2+\half(L_F M+F(M))^2\text{~~if~} u> M.
$$
Hence,
$$
\SK m(K)v_K F(u_K)\<\|v\|^2_{\LO}+\half L^2_F\|u\|^2_{\LO}+
\half C_F
$$
for all $v,u\in H_{\M}$, where $\stl{C_F=\half m(\Om)(\max_{0\<u\<M}F^2(u)+(L_F M+F(M))^2)}$. We obtain
$$
\bl
|A_R|&\stl \< \int^{T-\tau}_{0}\sum_{n\in n(t,\tau)}\dt(\|\uhk(\cdot,\lc t/\dt\rc)\|^2_{\LO}+\|\uhk(\cdot,{\lc (t+\tau)/\dt\rc})\|^2_{\LO})dt\\
     &\stl +\int^{T-\tau}_{0}\sum_{n\in n(t,\tau)}\dt ( L^2_F\|\uhk(\cdot,t_n)\|^2_{\LO}+C_F)dt.
\el
$$
One more time it follows from the estimate (\ref{ineq:energy}) and the Lemmas \ref{lemma:EGH1} and \ref{lemma:EGH2}
that
$$
|A_R|\<\tau \sum_{n=1}^{N}\dt (C\|\uhk(\cdot,t_n)\|^2_{\LO}+C_F)\<\tau C
$$
In the same way we proceed for the term $|A_S|$, one has that
$$
A_S=\int^{T-\tau}_{0}\sum_{n\in n(t,\tau)}\SK \int^{t_{n}}_{t_{n-1}}\IK (u^{\lc (t+\tau)/\dt\rc}_K-u^{\lc t/\dt\rc}_K) q(\x,s)\dx ds dt
$$
and
$$
\bl
|A_S|&\stl\<\int^{T-\tau}_{0}\sum_{n\in n(t,\tau)}\SK \int^{t_{n}}_{t_{n-1}}\IK \half((u^{\lc t/\dt\rc}_K)^2+(u^{\lc (t+\tau)/\dt\rc}_K)^2)+ (q(\x,s))^2\dx ds dt\\
&\stl =\int^{T-\tau}_{0}\sum_{n\in n(t,\tau)}\half\dt(\|\uhk(\cdot,{\lc t/\dt\rc})\|^2_{\LO}+\|\uhk(\cdot,\lc (t+\tau)/\dt\rc)\|^2_{\LO})dt\\
&\stl +\int^{T-\tau}_{0}\sum_{n\in n(t,\tau)}\SK \int^{t_{n}}_{t_{n-1}}\IK (q(\x,s))^2\dx ds dt.
\el
$$
In view of Lemmas \ref{lemma:EGH2} and \ref{lemma:EGH1} we obtain
$$
|A_S|\<\tau (\|\uhk\|^2_{\LQT}+\|q\|^2_{\LQT}).
$$
Finally we use an a priori estimate (\ref{ineq:energy}) and the hypothesis ($\hyp_6$) to conclude the proof.

\section{Estimate on space translates}
In this section we prove an estimate on the $L^2-$norm of differences of space
translates of the discrete solution.
We state without proof two results from
\cite{Eym_Gal_Her_2}, which are useful in our study.

\begin{Lem}\label{lemma:LqW1p}
Let $d\geqslant 1$, $1\< p<\infty$ and $\Om$ be an open
bounded connected subset of $\real^d$. Let $\D $ be a
mesh of $\Om$ in the sense of Definition \ref{Mesh}. Let
$\eta>0$ be such that $\eta\<
d_{K,\s}/d_{L,\s}\< 1/\eta$ for all
$\s\in\M_{\s}=\{K,L\}$. Then, there exists $q>p$
only depending on $p$ and there exist a positive constant $C$,
only depending on $d$, $\Om$, $p$ and $\eta$ such that:
\be\label{ineq:LqW1p}
\|u\|_{L^{q}(\Om)}\< C\|u\|_{1,p,\M}
\ee
for all $ u\in H_{\M}(\Om)$. We recall that $H_\M(\Om)\subset
\LO$ is the set of piecewise constant functions on the
control volumes of the mesh.
\end{Lem}
\begin{Lem}\label{lemma:L1est}
Let $d\geqslant 1$ and $\Om$ be a polyhedral open bounded connected
subset of $\real^d$. Let
$\D =(\M,\E,\mathcal{P})$ be a
discretization of $\Om$ in the sense of Definition \ref{Mesh}
and let $u\in H_{\M}(\Om)$. Then, with notation of Definition \ref{def:W1p}:
\be\label{ineq:L1est}
\|u(\cdot+\bold{y})-u\|_{L^1(\real^d)}\<|\bold{y}|\sqrt{d}\|u\|_{1,1,\M},
\ee
where $u$ is defined on the whole $\real^d$, taking $u=0$ outside $\Om$.
\end{Lem}
Next we show that a similar inequality holds in every $L^p$-norm.
\begin{Lem}\label{lemma:Lpest}
Let $d\geqslant 1$, $1\< p<\infty$ and $\Om$ be an open
bounded connected subset of $\real^d$ and $T>0$. Let
$\D $ be a discretization of $\Om$ in the sense of
Definition \ref{Mesh}. Let $\eta>0$ such that
$\eta\< d_{K,\s}/d_{L,\s}\< 1/\eta$ for all
$\s\in\M_{\s}=\{K,L\}$. There exist $C>0$
and $\rho>0$, only depending on $d$, $p$, $\Om$ and $\eta$ such
that
$$
\|u(\cdot+\bold{y})-u\|_{L^p(\real^d)}\<
C|\bold{y}|^{\rho}\|u\|_{1,p,\M},
$$
where $u$ is defined on $\real^d$, taking $u=0$ outside $\Om$.
\end{Lem}
\emph{Proof.}
In view of Lemma \ref{lemma:LqW1p}, there
exist $q>p$ and a positive constant $C$ such that
\be\label{loc:sp}
\|u\|_{L^q(\real^d)}\< C \|u\|_{1,p,\M}.
\ee
We apply the Interpolation Inequality [\cite{Adam_Four}, Theorem 2.11, p.27]
\be\label{loc:sp1}
\|u(\cdot+\bold{y})-u\|_{L^p(\real^d)}\<\|u(\cdot+\bold{y})-u\|^{\rho}_{L^1(\real^d)}\|u(\cdot+\bold{y})-u\|^{1-\rho}_{L^q(\real^d)},
\ee
where
$$
\rho=\frac{1}{p}\cdot\frac{q-p}{q-1}.
$$
Moreover (\ref{loc:sp}) implies that
$$
\|u(\cdot+\bold{y})-u\|_{L^q(\real^d)}\<2\|u\|_{L^q(\real^d)}\<C\|u\|_{1,p,\M},
$$
so that by (\ref{ineq:L1est}) and (\ref{loc:sp1}) implies that
$$
\|u(\cdot+\bold{y})-u\|_{L^p(\real^d)}\< C
|\bold{y}|^{\rho}(\|u\|_{1,1,\M})^{\rho}(\|u\|_{1,p,\M})^{1-\rho}.
$$
Applying H\"{o}lder inequality we obtain that
$$
\|u\|_{1,1,\M}\<C\|u\|_{1,p,\M}
$$
for some positive constant $C$. Then
$$
\|u(\cdot+\bold{y})-u\|_{L^p(\real^d)}\< C
|\bold{y}|^{\rho} \|u\|_{1,p,\M}.
$$
\begin{Th}\label{th:compactness}
Let $\{\D_h\}$ be a family of discretizations in the sense of
Definition \ref{Mesh} and let $\theta$ be a positive constant such that
$\theta_{\D}\<\theta$ for all $\D\in\D_h$. Let
$\{\uhk\}$ be a family of approximate solutions corresponding
to $\D_h$ and $\dt=T/N$ for some $N\in\intpos$. Then $\{\uhk\}$ is relatively compact in $\LQT$.
\end{Th}
\emph{Proof.} To begin with, we extend $\uhk$ by zero outside of $Q_T$. Applying the Lemma \ref{lemma:Lpest} with $p=2$
yields
$$
\|\uhk(\cdot+\bold{y},t)-\uhk(\cdot,t)\|_{L^2(\real^d)}\<
C |\bold{y}|^{\rho}\|\uhk(\cdot,t)\|_{1,2,\M}
$$
for some positive constants $\rho>0$ and $C>0$. Integrating on $(0,T)$ we
obtain
$$
\|\uhk(\cdot+\bold{y},\cdot)-\uhk\|^2_{L^2(\real^d\times (0,T))}\< C
|\bold{y}|^{2\rho}\Sn \dt\|\uhk(\cdot,t_n)\|^2_{1,2,\M}.
$$
Then in view of the lemmas \ref{lemma:w12bound},
\ref{lemma:equivalence} and the estimate (\ref{ineq:energy}) we
obtain the bound
$$
\|\uhk(\cdot+\bold{y})-\uhk\|_{L^2(\real^d\times (0,T))}\< C
|\bold{y}|^{\rho},
$$
which, combined with (\ref{ineq:TTrans}) gives
$$
\|\uhk(\cdot+\bold{y},\cdot+\tau)-\uhk\|_{L^2(\real^d\times (0,T))}
$$
$$
\<\|\uhk(\cdot+\bold{y},\cdot+\tau)-\uhk(\cdot+\bold{y},\cdot)\|_{L^2(\real^d\times (0,T))}+\|\uhk(\cdot+\bold{y},\cdot)-\uhk\|_{L^2(\real^d\times (0,T))}\\
$$
$$
\< C(\sqrt{\tau}+|\bold{y}|^\rho).
$$
Then the Fr\'{e}chet-Kolmogorov Compactness Theorem implies that
the family $\{\uhk\}$ is relatively compact in $L^2(\real^d\times (0,T))$ and thus in $\LQT$.

\section{Convergence}
\begin{Th}\label{th:weak_conv_grad}
Let $\{\D_h\}$ be a family of discretizations in the sense of
Definition \ref{Mesh} and let $\theta$ be a positive constant such that
$\theta_{\D}\<\theta$ for all $\D\in\D_h$. Let
$\{\uhk\}$ be a family of approximate solutions corresponding
to $\D _h$ and $\dt=T/N$ for some $N\in\intpos$. Then there exists a subsequence of $\{\uhk\}$, which we denote again by $\{\uhk\}$, such that $\uhk\to u$ strongly in $\LQT$ as $\hD,\dt\to 0$, where $u$ is a weak solution of Problem $\Pb$. Moreover $u\in L^2(0,T;H^1_0(\Om))$ and $\grad_{\D,\dt}\uhk$ weakly converge in $\LQT^d$ to $\grad u$. In the case that $F$ is nondecreasing, the whole sequence $\{\uhk\}$ converges to the unique weak solution $u$ of Problem $\Pb$.
\end{Th}

\emph{Proof.} By Theorem \ref{th:compactness} there exist a subsequence of
$\{\uhk\}$ that we still denote by $\{\uhk\}$ and a function $u\in
\LQT$ such that $\uhk\to u$ strongly in $\LQT$
as $\hD,\dt\to 0$ (and also in $L^2(\real^d\times (0,T))$ taking
$\uhk=0$ outside of $\Om\times (0,T)$). In view of (\ref{ineq:energy}) there exists
a function $\bold{G}\in \LQT^d$ such that
$\grad_{\D,\dt}\uhk$ weakly converge in $\LQT^d$ to $\bold{G}$ along a subsequence
as $\hD,\dt\to 0$. In order to show that $\bold{G}=\grad u$
we consider an arbitrary vector function $\bold{\w}\in C([0,T];C^{\infty}_c(\real^d))$
and the term $T^1_{\bold{G}}$ defined by
$$
T^1_{\bold{G}}=\IT \int_{\real^d}\grad_{\D,\dt}\uhk(\x,t)\cdot
\bold{\w}(\x,t)\dx dt.
$$
Using the Definition \ref{def:solution} and (\ref{D_K,s}) we
obtain that $T^1_{G}=T^2_{G}+T^3_{G}$, with
$$
T^2_{\bold{G}}=\SnKs \dt m(\s)(u^n_{\s}-u^n_{K})\n_{K,\s}\cdot\bold{\w}^n_K
$$
and
$$
T^3_{\bold{G}}=\SnKs R_{K,\s}(u^n)\n_{K,\s}\cdot\Itn \int_{D_{K,\s}}\bold{\w}(\x,t)\dx dt,
$$
where
$\stl \bold{\w}^n_K=\frac{1}{\dt m(K)}\ItnK  \bold{\w}(\x,t)\dx dt$.
We compare $T^2_{\bold{G}}$ with $T^4_{\bold{G}}$ defined by
$$
T^4_{\bold{G}}=\SnKs \dt m(\s)(u^n_{\s}-u^n_{K})\n_{K,\s}\cdot\bold{\w}^n_{\s},
$$
where
$\stl{\bold{\w}^n_{\s}=\frac{1}{\dt m(\s)}\Itn \int_{\s}\bold{\w}(\x,t)~d\gamma dt}$.
One can see that
$$
(T^2_{\bold{G}}-T^4_{\bold{G}})^2\<\SnKs \frac{\dt m(\s)}{d_{K,\s}}(u^n_{\s}-u^n_{K})^2\SnKs \dt m(\s)d_{K,\s}|\bold{\w}^n_{K}-\bold{\w}^n_{\s}|^2,
$$
which leads to $T^2_{\bold{G}}\to T^4_{\bold{G}}$ as
$\hD\to 0$. Then,
$$
T^4_{\bold{G}}=-\SnKs \dt m(\s)u^n_{K}\n_{K,\s}\cdot\bold{\w}^n_{\s}=-\IT \int_{\real^d}\uhk(\x,t)~\dv \bold{\w}(\x,t)\dx dt.
$$
and we conclude that
$$
\stl{\lim_{\hD,\dt\to
0}T^4_{\bold{G}}=-\IT \int_{\real^d}u(x,t)\dv ~\bold{\w}(\x,t)\dx dt}.
$$
Next, we show that
$\stl{\lim_{\hD,\dt\to 0}T^3_{\bold{G}}=0}$. Thanks to (\ref{eq:E(RKs)}) we
have that
$$
T^3_{\bold{G}}=\SnKs R_{K,\s}(u^n)\n_{K,\s}\cdot\Itn \int_{D_{K,\s}}(\bold{\w}(\x,t)-\bold{\w}^n_K)\dx dt.
$$
Since $\bold{\w}$ is a regular function, there exist a positive
constant $C=C(\bold{\w})$ such that
$$
\stl{|\Itn \int_{D_{K,\s}}(\bold{\w}(\x,t)-\bold{\w}^n_K)\dx dt}|\<
C \dt \frac{m(\s)d_{K,\s}}{d}(\hD+\dt).
$$
On the other hand by (\ref{R_K,s}) and in view of regularity of
the mesh, we have that
$$
(R_{K,\s}u^n)^2\<
2d((\frac{u^n_{\s}-u^n_{K}}{d_{K,\s}})^2+|\grad_{K}u^n|^2|\frac{\x_{\s}-\x_{K}}{d_{K,\s}}|^2)\<
2d((\frac{u^n_{\s}-u^n_{K}}{d_{K,\s}})^2+\theta^2|\grad_{K}u^n|^2).
$$
Applying the \CaSz inequality, we get
$$
\lim_{\hD,\dt\to 0}T^3_{\bold{G}}=0
$$
This implies that the function $\bold{G}\in L^2(\real^d\times (0,T))^d$ is a.e.
equal to $\grad u$ in $\real^d\times (0,T)$. Since $u=0$ outside of $\Om$, it follows that $u\in L^2(0,T;H^1_0(\Om))$.\\
Next we show that $u$ is a weak solution of the problem $\Pb$. For
this purpose, we introduce the function space
$$
\ba{lcr}
\Phi=\{\ph\in C^{2,1}(\overline{\Om}\times
[0,T]),&\ph=0$ on $\partial\Om\times
[0,T],&\ph(\cdot,T)=0\}.
\ea
$$
Taking an arbitrary $\ph\in \Phi$, we define the sequence of elements of $X_{\D ,0}$
$$
\ph^n=P_{\D}\ph(\cdot,t_n) \tfa n\in\{1,\ldots,N\}
$$
which implies $\ph^n_K=\ph(\x_K,t_n)$ and $\ph^n_{\s}=\ph(\x_{\s},t_n)$. Next setting
$$
v^n=\ph^{n-1} \tfa n\in\{1,\ldots,N\}
$$
in (\ref{eq:var}), we obtain, also in view of (\ref{bil_F:var}) and (\ref{bil_T:var}) that
$$
T_T+T_D+T_C+T_R=T_S,
$$
where
$$
T_T=\SnK m(K)(\beta(u^n_K)-\beta(u^{n-1}_K))\ph^{n-1}_K,
$$
$$
T_D=\Sn \dt\SK \Ss (\ph^{n-1}_K-\ph^{n-1}_{\s})F_{K,\s}(u^n),
$$
$$
T_C=\Sn \dt\SK \Ss (\ph^{n-1}_K-\ph^{n-1}_{\s})
V_{K,\s}\ubar,
$$
$$
T_R=\Sn \dt\SK m(K)\ph^{n-1}_K F(u^n_{K})
$$
and
$$
T_S=\Sn \dt\SK m(K)\ph^{n-1}_K q^n_{K}.
$$
We successively search for the limit of each of these terms as $\hD$ and $k$ tend to zero.
\subsection{Time evolution term} Using discrete integration by
parts and the fact that $\ph(\x,T)=0$ we obtain
$$
T_T=-\SnK m(K)(\ph^n_K-\ph^{n-1}_K)\beta(u^n_K)-\SK m(K)\ph^0_K
\beta(u^0_K).
$$
First we show that
$$
\SK m(K)\beta(u^0_K) \ph^0_K \to
\IO \beta(u_0(\x))\ph(\x,0)\dx.
$$
For this purpose we define
$$
T^0_T=\SK m(K)\beta(u^0_K) \ph^0_K-\IO \beta(u_0(\x))\ph(\x,0)\dx.
$$
Next we subtract $\IO \beta(u^0_K)\ph(\x,0)$ from each term to deduce that,
\be\label{loc:convergence:T0}
T^0_T=\SK \IK \beta(u^0_K)(\ph^0_K-\ph(\x,0))\dx-\SK \IK (\beta(u_0(\x))-\beta(u^0_K))\ph(\x,0)\dx.
\ee
In view of the regularity of the test function $\ph\in C^{2,1}(\overline{\Om}\times[0,T])$ we have that
$$
|\ph^0_K-\ph(\x,0)|\< C\hD \tfa \x\in K
$$
and
$$
|\ph(\x,0)|\< C.
$$
We also remark that by (\ref{eq:disc:ini}) one has that $|u^0_K|\<\|u_0\|_{L^\infty(\Om)}$ and moreover the monotonicity hypothesis ($\hyp_1$) implies that $|\beta(u^0_K)|\< \beta(\|u_0\|_{L^\infty(\Om)})$ for all $K\in\M$; consequently the first term on the \rhs of (\ref{loc:convergence:T0}) tends to zero as $\hD\to 0$ and the second term can be estimated by
$$
C\SK \IK |\beta(u_0(\x))-\beta(u^0_K)|\dx = C\IO |\beta(u_0(\x))-\beta(\uhk(\x,0))|\dx.
$$
By the discrete initial condition $(\ref{eq:disc:ini})$ and the Definition $\ref{def:solution}$ one has
$$
\IO |u_0(\x)-\uhk(\x,0)|\dx\to 0\text{~as~}\hD\to 0,
$$
or in other words $\uhk(0)$ converges strongly to $u_0$ in $L^1(\Om)$ as $\hD\to 0$.
Hence a subsequence of $\{\uhk(\x,0)\}$, which we still denote by $\{\uhk(\x,0)\}$ converges to $u_0(\x)$ for a.e. $x\in\Om$ and also $\beta(\uhk(\x,0))\to\beta(u_0(\x))$ for a.e. $x\in\Om$. Since $\beta(\uhk(\x,0))\<\|\beta(u_0(\x))\|_{L^{\infty}(\Om)}$ the Lebesgue dominated convergence theorem implies
$$
\IO |\beta(u_0(\x))-\beta(\uhk(\x,0))|\dx\to 0\text{~as~}\hD\to 0.
$$
Thus $T^0_T\to 0\text{~as~}\hD\to 0$. Next we prove that
\be
\label{loc:convergence:time_1}
\SnK m(K)(\ph^n_K-\ph^{n-1}_K)\beta(u^n_K)\to
\ITO \beta(u(\x,t))\ph_{t}(\x,t)\dx dt
\ee
as $\hD$ and $\dt\to 0$. We define
$$
T^1_T=\SnK m(K)(\ph^n_K-\ph^{n-1}_K)\beta(u^n_K)-
\ITO \beta(u(\x,t))\ph_{t}(\x,t)\dx dt,
$$
and we add and subtract
$\stl{\ItnK  \beta(u^n_K)\ph_{t}(\x,t)\dx dt}$
in each term to obtain
\be\label{T^1_T}
\bl
T^1_T
&\stl=\SnK m(K)\beta(u^n_K)\Itn (\ph_{t}(\x_K,t)-\ph_{t}(\x,t))\dx dt\\
&\stl+\ITO (\beta(\uhk(\x,t))-\beta(u(\x,t)))\ph_{t}(\x,t)\dx dt.
\el
\ee
We have that for all $x\in K$ and all $K\in\M$ it holds
$$
|\ph_{t}(\x_K,t)-\ph_{t}(\x,t)|\< C(\hD)
$$
where $C(\hD)\to 0$ as $\hD\to 0$. The absolute value of the first term on the \rhs of (\ref{T^1_T}) is bounded by
$$
\bl
\stl C(\hD)\Sn \SK \dt m(K)|\beta(u^n_K)|
&\stl\<C(\hD) (T m(\Om))^{1/2}(\Sn \SK \dt m(K)(\beta(u^n_K))^2)^{1/2}\\
&\stl\< C(\hD) T m(\Om)^{1/2}\|\beta(\uhk)\|_{\LinfO},
\el
$$
which tends to zero as $\hD\to 0$ in view of the a priori estimate (\ref{ineq:b(u)}). Further, since $|\ph_{t}(\x,t)|\< C_{\ph}$, we can
estimate the absolute value of the second term in (\ref{T^1_T}) by
$$
\bl
\stl C_{\ph}\ITO |\beta(\uhk(\x,t))-\beta(u(\x,t))|\dx dt
&\stl\<C_{\ph}\ITO |\widetilde{\beta}_1(\uhk(\x,t))-\widetilde{\beta}_1(u(\x,t))|\dx dt\\
&\stl+C_{\ph}\ITO |\widetilde{\beta}_2(\uhk(\x,t))-\widetilde{\beta}_2(u(\x,t))|\dx dt,
\el
$$
where $\widetilde{\beta}_1$ and $\widetilde{\beta}_2$ are given by (\ref{def:splitting1})-(\ref{def:splitting3}).
Since $\uhk\to u$ strongly in $\LQT$, a subsequence of $\{\uhk\}$, which we still denote by $\{\uhk\}$ converges to $u$ a.e. in $\Om$.
The first term on \rhs of the expression above converges to zero by the Lebesgue dominated convergence theorem. The convergence to zero of the second term can be deduced from the \Lip continuity of $\widetilde{\beta}_2$ and the strong convergence of $\uhk$ to $u$ in $\LQT$.
\subsection{Convection term}
Next, we show that
\be\label{loc:convergence:conv_3}
E = \Sn \dt\SK \Ss (\ph^{n-1}_K-\ph^{n-1}_{\s})
V_{K,\s}\ubar + \ITO u(\x,t)\bold{V(\x)}\cdot\grad\ph(\x,t)\dx dt \to 0
\ee
as $\hD,\dt$ tend to zero. As it was done in in the proof of the Lemma \ref{lemma:bils} we write the left-hand side part of (\ref{loc:convergence:conv_3}) as
\be
\ba{c}\label{loc:convergence:conv_2}
\stl  \Sn \dt\SK \Ss (\ph^{n-1}_K-\ph^{n-1}_{\s})V_{K,\s}\ubar
\stl =\Sn \dt\SK \Ss V_{K,\s} (\ph^{n-1}_K-\ph^{n-1}_{\s})u^n_K\\
\stl-\Sn \dt\SK \sum_{\s\in\E_K, V_{K,\s}\<0}V_{K,\s}(\ph^{n-1}_K-\ph^{n-1}_{\s})(u^n_K-u^n_{\s}),
\ea
\ee
and the following estimate holds
$$
|\Sn \dt\SK \sum_{\s\in\E_K, V_{K,\s}\<0}V_{K,\s}(\ph^{n-1}_K-\ph^{n-1}_{\s})
(u^n_K-u^n_{\s})|\<\hD\cdot\Vmax \Sn \dt|\ph^{n-1}|_X |u^n|_X
$$
$$
\<C\hD\cdot\Vmax \Sn \dt\|\grad_{\D }\ph^{n-1}\|_{\LO} \|\grad_{D}u^n\|_{\LO}
$$
$$
\<C\hD\cdot\Vmax \Sn \dt\|\grad_{\D }\ph^{n-1}\|^2_{\LO} +C\hD\cdot\Vmax \|\grad_{\D,\dt} \uhk\|^2_{\LQT}.
$$
The second term in the \rhs of the expression above is bounded be because of the a priori estimate (\ref{ineq:energy}) and the first term can be controlled via the consistency of the discrete gradient given by Lemma \ref{lemma:DiscGradLInf} and the regularity of $\ph$; indeed
$$
\Sn\Itn\IO |\grad_{\D }\ph^{n-1}(\x)|^2\dx dt
\<3\Sn\Itn\IO |\grad_{\D }\ph^{n-1}(\x)-\grad \ph(\x,t_{n-1})|^2\dx dt
$$
$$
+3\Sn\Itn\IO |\grad\ph(\x,t_{n-1})-\grad \ph(\x,t)|^2
\dx dt
+3\ITO
|\grad \ph(\x,t)|^2
\dx dt
$$
$$
\<(C_{\ph}(\dt)+C\hD^2) T m(\Om)+3\|\grad\ph\|^2_{\LQT}\< C,
$$
where $C_{\ph}(\dt)$ tends to zero as $\dt\to 0$. Thus the second term in the right hand side of (\ref{loc:convergence:conv_2}) tends to zero as $\hD,\dt\to 0$. Let us define $E_1$ and $E_2$
$$
E_1 = \Sn \dt \SK \Ss V_{K,\s}(\ph^{n-1}_K-\ph^{n-1}_{\s})u^n_K+
\Sn \dt\SK u^n_K\IK \grad\ph(\x,t_{n-1})\cdot\bold{V(\x)}\dx,
$$
$$
E_2 = \Sn \dt\SK u^n_K\IK \grad\ph(\x,t_{n-1})\cdot\bold{V(\x)}\dx-
\ITO u(\x,t)\grad\ph(\x,t)\cdot\bold{V(\x)}\dx dt.
$$
so that also in view of (\ref{loc:convergence:conv_2}) $\lim_{\hD,\dt\to 0}E = \lim_{\hD,\dt\to 0} E_1 - \lim_{\hD,\dt\to 0} E_2$. We will successively establish that $E_1$ and $E_2$ converges to zero as $\hD,\dt\to 0$.
To begin with let us remark that integrating by parts in the expression of $E_1$ yields $E_1 = E_{11}-E_{12}$, where
$$
E_{11} = \Sn \dt
\SK \Ss V_{K,\s}\ph^{n-1}_K u^n_K
-\Sn \dt\SK u^n_K\IK \ph(\x,t_{n-1})\dv\V(\x) \dx
$$
and where
$$
E_{12} = \stl\Sn \dt\SK \Ss V_{K,\s}\ph^{n-1}_{\s} u^n_K-\Sn \dt\SK \Ss u^n_K\int_{\s}\ph(\x,t_{n-1})\V(\x)\cdot \n_{K,\s} d\gamma.
$$

We first prove that $\stl \lim_{\hD,\dt\to 0}E_{11} = 0$.
$$
E_{11} = \Sn \dt
\SK u^n_K\IK (\ph^{n-1}_K - \ph(\x,t_{n-1}))\dv\V(\x))\dx.
$$
in view of regularity of the function $\ph$ we obtain
$$
|E_{11}| \< C_{\ph} \hD \Sn \dt
\SK |u^n_K|\IK |\dv \V(\x)|\dx\<C_{\ph} h \int_{Q_T}|\uhk(\x,t) \dv \V(\x)|\dx
$$
Finally applying the \CaSz yields
$$
|E_{11}| \<C_{\ph} h \|\uhk\|_{\LQT} \|\dv \V\|_{\LQT}
$$
so that $|E_{11}|\to 0$ as $\hD\to 0$. Next we consider the term $E_{12}$, which can be written as
$$
\stl E_{12}= \Sn \dt\SK \Ss u^n_K\int_{\s}(\ph^{n-1}_{\s}-\ph(\x,t_{n-1}))\V(\x)\cdot \n_{K,\s} d\gamma.
$$
In order to show that $E_{12}\to 0$ as $\hD,\dt\to 0$ we first remark that since $\n_{K,\s}=-\n_{L,\s}$ for any pair of neighbor volumes $K, L$, and in view of the boundary condition on $\ph$ one has that
$$
\Sn \dt
\SK \Ss V_{K,\s}\ph^{n-1}_{\s} u^n_{\s}=0
$$
and also
$$
\Sn \dt\SK \Ss u^n_{\s}\int_{\s}\ph(\x,t_{n-1})\V(\x)\cdot \n_{K,\s} d\gamma=0.
$$
Hence, the term $E_{12}$ can be written as
$$
E_{12} = \Sn \dt\SK \Ss (u^n_K-u^n_{\s})\int_{\s}(\ph^{n-1}_{\s}-\ph(\x,t_{n-1}))\V(\x)\cdot \n_{K,\s} d\gamma.
$$
Therefore, in view of the regularity of $\ph$ and $\V$ we have that
$$
|E_{12}| \< C\max_{\x\in\Om}\V(\x) \cdot h\Sn \dt\SK \Ss m(\s)|u^n_K-u^n_{\s}|.
$$
Applying \CaSz we obtain
$$
|E_{12}| \< C d^{\half}\max_{\x\in\Om}\V(\x) \cdot \hD(\Sn \SK \Ss  \dt m(\s)\frac{(u^n_K-u^n_{\s})^2}{d_{K,\s}})^{\half}\cdot (\Sn \SK \Ss  \dt \frac{m(\s)d_{K,\s}}{d})^{\half}.
$$
In view of Lemma \ref{lemma:equivalence} we obtain
$$
|E_{12}|\< C d\max_{\x\in\Om}\V(\x) m(\Om)^{\half}T^{\half}\cdot \hD\|\grad_{\D,\dt} \uhk\|_{L(Q_T)},
$$
so that in view of the a priori estimate (\ref{ineq:energy})
one has that $|E_{12}|\to 0$ as $\hD\to 0$, so that $E_1=E_{11}-E_{12}\to 0$ as $\hD,\dt\to 0$.
It remains to prove that $E_2$ converges to zero.
\par
Adding and subtracting
$\stl{\ItnK  u^n_K\grad\ph(\x,t)\cdot\bold{V(\x)}\dx dt}$
from each term of $E_2$, yields
$$
\bl
E_2
&\stl=\SnK \ItnK  u^n_K(\grad\ph(\x,t_{n-1})-\grad\ph(\x,t))\cdot\bold{V(\x)}\dx dt\\
&\stl-\SnK \ItnK  (u(\x,t)-\uhk(\x,t))\grad\ph(\x,t)\cdot\bold{V(\x)}\dx dt.
\el
$$
Finally, in view of the regularity of $\ph$, the a priori estimate (\ref{ineq:energy}), and to the fact
that $\uhk\to u$ strongly in $\LQT$ we conclude
that $|E_2|$ tends to zero as $\hD,\dt\to 0$.
\subsection{Diffusion term}
We show below that
$$
T^1_D=\Sn \dt\SK \Ss (\ph^{n-1}_K-\ph^{n-1}_{\s})F_{K,\s}(u^n)
-\IT\IO\grad \ph(\x,t)\cdot
\Lam(\x)\grad u(\x,t)\dx dt
$$
tends to zero as $\hD,\dt\to 0$. In view of (\ref{bil_Flux})
one has that
$$
T^1_D=\Sn \Itn \IO (\grad_{\D }\ph^{n-1}
\cdot \Lam(\x)\grad_{\D } u^n-\grad
\ph(\x,t)\cdot \Lam(\x)\grad u(\x,t))\dx dt.
$$
Adding and  subtracting the term $\stl{\IO \grad
\ph(\x,t)\cdot\Lam(\x)\grad_{\D }u^n}\dx$
we set $T^1_D$ in the form $T^1_D=T^2_D+T^3_D$ with
$$
T^2_D=\Sn \Itn \IO (\grad_{\D }\ph^{n-1}-\grad\ph(\x,t))
\cdot\Lam(\x)\grad_{\D }u^n\dx dt
$$
and
$$
T^3_D=\Sn \Itn \IO \grad\ph(\x,t)
\cdot\Lam(\x)(\grad_{\D }u^n-\grad
u(\x,t))\dx dt.
$$
The term $T^3_D$ tends to zero as $\hD,\dt\to 0$, since $\grad_{\D,\dt} \uhk$ tends to $\grad u$ weakly in $\LQT$. On the other hand the term $T^2_D$
can be written in the form $T^2_D=T^4_D+T^5_D$ with
$$
T^4_D=\Sn \Itn \IO (\grad_{\D }\ph^{n-1}-\grad\ph(\x,t_{n-1}))
\cdot\Lam(\x)\grad_{\D }u^n\dx dt
$$
and
$$
T^5_D=\Sn \Itn \IO \grad(\ph(\x,t_{n-1})-\ph(\x,t))
\cdot\Lam(\x)\grad_{\D }u^n\dx dt.
$$
It follows from (\ref{ineq:energy}), Lemma \ref{lemma:DiscGradLInf} and the
regularity of $\ph$ that $T^4_D$ and $T^5_D$ tends to
zero as $\hD,\dt\to 0$ and so do $T^2_D$ and $T^1_D$.
\subsection{Reaction term}
Let us show that
$$
T_R\to\ITO F(u(\x,t))\ph(\x,t)\dx dt
$$
as $\hD$ and $k$ tend to zero. For this purpose, we introduce
$$
T^1_R=\SnK \ItnK  (\ph^{n-1}_K-\ph(\x,t))F(u^n_{K})\dx dt,
$$
$$
T^2_R=\SnK \ItnK  \ph(\x,t)( F(u^n_{K}) - F(u(\x,t)))\dx dt.
$$
We obtain the convergence result similarly as for the time evolution term; more precisely we split the reaction $F$ into a bounded and a \Lip continuous parts by setting
\begin{center}

$
F_1(s)=
\left\{\ba{cc}
F(s)& 0\<s\< M\\
0&~\text{otherwise},
\ea\right.
$
$
F_2(s)=
\left\{\ba{cc}
0&0\<s\< M\\
F(s)&~\text{otherwise},
\ea\right.
$
\end{center}
and
$$
y(s)=
\left\{\ba{cc}
\stl{\frac{F(M)}{M}s} & 0\<s\< M\\
\stl{0}&~\text{otherwise}.
\ea\right.
$$
We, then, define $\widetilde{F}_1=F_1-y$ and $\widetilde{F}_2=F_2+y$ which are both continuous; moreover $|\widetilde{F}_1|$ is bounded by $C_{\widetilde{F}}=\stl\max_{0\<s\<M}|F(s)|+F(M)$, and $\widetilde{F}_2$ is \Lip continuous with \Lip constant $L_{\widetilde{F}}=\stl{\max(L_{F},F(M)/M)}$.
In view of the regularity of $\ph$ one has
$$
|\ph(\x,t)-\ph^{n-1}_K|\<C(\hD+\dt)\tfa x\in K, t\in(t_{n-1},t_n],
$$
so that
$$
|T^1_R|\<C(\hD+\dt)\ITO |\widetilde{F}_1(\uhk(\x,t))+\widetilde{F}_2(\uhk(\x,t))|\dx dt
$$
$$
\<C(\hD+\dt)(C_{\widetilde{F}}T m(\Om)+L_{\widetilde{F}}T^{\half} m(\Om)^{\half}\|\uhk\|_{\LQT}),
$$
which by the a priori estimate (\ref{ineq:energy}) implies that $|T^1_R|\to 0$ as $\hD,\dt \to 0$. Since $\ph$ is bounded we can estimate the second term as
$$
|T^2_R|\<C \ITO |F(\uhk(\x,t)) - F(u(\x,t))|\dx dt
$$
$$
\<\ITO C|\widetilde{F}_1(\uhk(\x,t))-\widetilde{F}_1(u(\x,t))|\dx dt+
\ITO C|\widetilde{F}_2(\uhk(\x,t)) - \widetilde{F}_2(u(\x,t))|\dx dt.
$$
The convergence is can be proved by applying similar arguments as for the time evolution term.
\subsection{Source term}
We deduce from the regularity of $\ph$ that
$$
T_S=\Sn \SK \ItnK  \ph(\x_K,t_{n-1})q(\x,t)\dx dt\to\ITO \ph(\x,t)q(\x,t)\dx dt
$$
as $\hD,\dt\to 0$.
\subsection{Convergence to a weak solution of Problem $\Pb$}
In view of Theorem \ref{th:weak_conv_grad} $\{\uhk\}$ strongly converges to $u$ in $L^2(Q_T)$, with $u \in L^2(0,T;H^1_0(\Om))$, and it follows from (\ref{ineq:b(u)}) that $\beta(u)\in\LinfO$. Moreover we deduce
from the density of the set $\Phi$ in the set $\{\ph\in L^2(0,T;H^1_0(\Om)), \ph_t\in L^{\infty}(Q_T),
\ph(\cdot,T)=0\}$ that $u$ is a weak solution of the continuous problem $\Pb$ in the sense of Definition \ref{def:weak solution}. In the case that $F$ is nondecreasing so that the solution of Problem $\Pb$ is unique (cf. Remark \ref{rq:unique}) we conclude that the whole family $\{\uhk\}$ converges to $u$.

\section{Numerical simulations}
In this section we present the results of numerical simulations. The purpose is to test our
scheme in the case of problems with a known analytical solution.\\
\subsection{Numerical Test I}
\begin{figure}[htb]
\begin{center}
\includegraphics[width=0.7\textwidth]{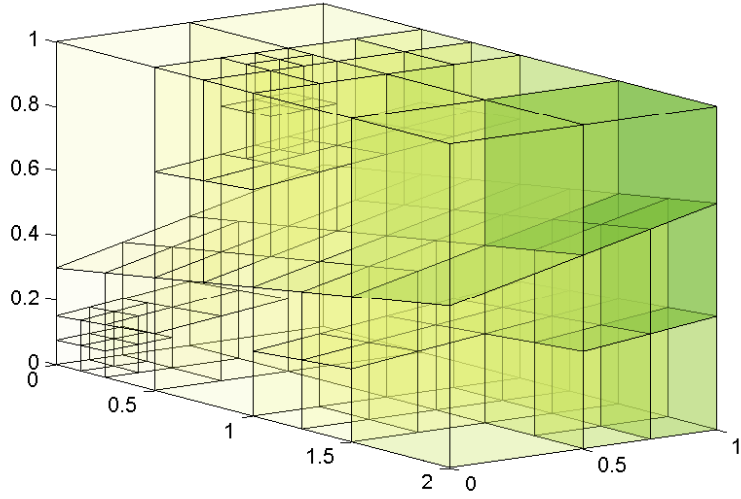}
\caption{Approximate solution on the nonmatching hexahedral mesh at $t=1$}\label{mesh}
\end{center}
\end{figure}

We consider the equation
$$
\ddt({(u+u^{\half})}) -\dv(\Lam(\x)\grad u)+\dv(\bold{V(\x)}u)+\half u^{\half} = 0
$$
in the $3$-dimensional space domain $\Om=(0,2)\times(0,1)\times(0,1)$. We define the discontinuous $\Lam$ and $\V$ fields as follows
\\
For all $x_1\<1$ we set
$$\Lam=\left(\ba{ccc}
1&0&0\\
0&1&0\\
0&0&1
\ea\right)~~\text{and}~~\V=(4,0,0);
$$
for all $x_1>1$ we set
$$\Lam=\left(\ba{ccc}
8&-5&-2\\
-5&20&-7\\
-2&-7&19
\ea\right)~~\text{and}~~\V=(4,7,7).
$$
The initial and the \Dir boundary conditions are given by the exact solution
$$
\stl{u(\x,t)=e^{x_1+x_2+x_3-t-3}}.
$$
We remark that the velocity field $\V$ and the total flux $\Lam(\x)\grad u +\V(\x)u$ have a continuous normal trace across the discontinuity $x=1$. We perform the simulations on $3$-dimensional hexahedral meshes with random refinement (see Figure \ref{mesh}), so that the mesh is nonmatching. In Table \ref{Tab} we present simulation results with various mesh sizes $\hD$ and time steps $k$; we denote by $Err$ the maximum relative error in $L^2$-norm, namely

$$
\stl Err=\max_{n\in\{1,\ldots,N\}}\frac{\|u_{h,t}(\cdot,t_n)-u(\cdot,t_n)\|_{\LO}}{\|u(\cdot,t_n)\|_{\LO}}.
$$

\begin{table}[htb]
\begin{center}
\begin{tabular}{|c|c|c|c||c|}
 N &h & \# of elements & \# of faces& Err\\
\hline
 50 & 0.75& 165  & 672 & 0.03575 \\
 100 & 0.375& 837 & 3324 & 0.01432 \\
 200 & 0.1875 & 3203 & 11550 & 0.00648  \\
 400 & 0.0938& 18533 & 60633 & 0.00305  \\
\end{tabular}
\end{center}
\caption{Number of time steps $N$, mesh diameter $\hD$, number of
elements, number of faces and the relative error for nonmatching hexahedral meshes}\label{Tab}
\end{table}
\subsection{Numerical test II}
We consider a degenerate parabolic equation which possesses a traveling wave solution, namely
$$
\ddt({(u^{\half})}) -\dv(\delta\grad u)+\dv((v,0,0)u) = 0
$$
in the domain
$$\Om=(0,1)^3 \ta T=1.$$
This equation admits the following $1$-dimensional exact solution
$$
u(x,y,t)=\stl{(1-e^{\frac{v}{2\delta}(x-vt-p)})^2} \tf x\< vt+p,
$$
$$
u(x,y,t)=0  \tf x>vt+p
$$
where $p, v, \delta$ are parameters still to be defined. We set $p=0.2$, $v=0.8$, and consider two values of $\delta$, namely $0.01, ~0.0001$. The initial state is given by the exact solution at the time $t = 0$ and we prescribe corresponding  \Dir boundary conditions on the sides $x=0$ and $x=1$. The null flux boundary condition is imposed on the remaining part of the boundary.

Since the scheme does not preserve the maximum principle, it is necessary to define the function $\beta(u)=\stl{u^{\half}}$ for negative values as well, which leads us to set $\beta(-u)=-\beta(u)$. Further one has to solve the system of nonlinear equations
\be\label{loc:test1}
\left\{\ba{rc}
\stl m(K)(\beta(u^n_{K})-\beta(u^{n-1}_{K}))
+k\Ss  {F_{K,\s}(u^n) }&\\
\stl +k\Ss  V_{K,\s}\ubar =\dt~m(K)q^n_K,& \tfa K\in\M,\\
(F_{K,\s}(u^n)+V_{K,\s}\ubar )+(F_{L,\s}(u^n)+V_{L,\s}\overline{u^n_{L,\s}})=0, &\tfa \s\in\E_{int}\\\\
u^n_{\s}=0, &\tfa \s\in\E_{ext}.
\ea\right.
\ee
Since $\beta'(0)=+\infty$ the Newton method can not be directly applied. In order to overcome this difficulty we introduce new discrete unknowns
$$
w^n=\beta(u^n),\text{~and thus~} u^n = \ph(w^n), \text{~where~} \ph = \beta^{-1}.
$$
In view of (\ref{upwind_Vu}) and (\ref{def:Flux}) the nonlinear system becomes
\be\label{nls}
\left\{\ba{rl}
\stl{m(K)(w^n_{K}-w^{n-1}_{K})}
+\stl \dt\sum_{\s,\s'\in\E_{K}}A^{\s \s'}(\ph(w^n_K)-\ph(w^n_{\s'})) &\\
+\stl \dt\Ss (V^+_{K,\s}\ph(w^n_K)+V^-_{K,\s}\ph(w^n_{\s}))=\dt~m(K)q^n_K,& \tfa K\in\M,\\

\stl\sum_{\s'\in\E_{K}}A^{\s \s'}(\ph(w^n_K)-\ph(w^n_{\s'}))+(V^+_{K,\s}\ph(w^n_K)+V^-_{K,\s}\ph(w^n_{\s}))&\\
+\stl\sum_{\s'\in\E_{L}}A^{\s \s'}(\ph(w^n_L)-\ph(w^n_{\s'}))+(V^+_{L,\s}\ph(w^n_L)+V^-_{L,\s}\ph(w^n_{\s}))=0, &\tfa \s\in\E_{int},\\

\ph(w^n_{\s})=0, &\tfa \s\in\E_{ext}.
\ea\right.
\ee
The system (\ref{nls}) depends on $(w^n_{\s})_{\s\in\E}$ only through the terms $(\ph(w^n_{\s}))_{\s\in\E}$. This lead us to choose the discrete unknowns
$$
w^n_K=\beta(u^n_K) \tfa K\in\M \ta u^n_{\s} \tfa \s\in\E,
$$
so that the system (\ref{nls}) takes the form
\be
\left\{\ba{rl}
\stl{m(K)(w^n_{K}-w^{n-1}_{K})}
+\stl \dt\sum_{\s,\s'\in\E_{K}}A^{\s \s'}(\ph(w^n_K)-u^n_{\s'})) &\\
+\stl \dt\Ss (V^+_{K,\s}\ph(w^n_K)+V^-_{K,\s}u^n_{\s})=\dt~m(K)q^n_K,& \tfa K\in\M,\\

\stl\sum_{\s'\in\E_{K}}A^{\s \s'}(\ph(w^n_K)-u^n_{\s'})+(V^+_{K,\s}u^n_K+V^-_{K,\s}u^n_{\s})&\\
+\stl\sum_{\s'\in\E_{L}}A^{\s \s'}(\ph(w^n_L)-u^n_{\s'})+(V^+_{L,\s}\ph(w^n_L)+V^-_{L,\s}u^n_{\s})=0, &\tfa \s\in\E_{int},\\

u^n_{\s}=0, &\tfa \s\in\E_{ext}.
\ea\right.
\ee

\begin{figure}[!htb]
\includegraphics[width=\textwidth]{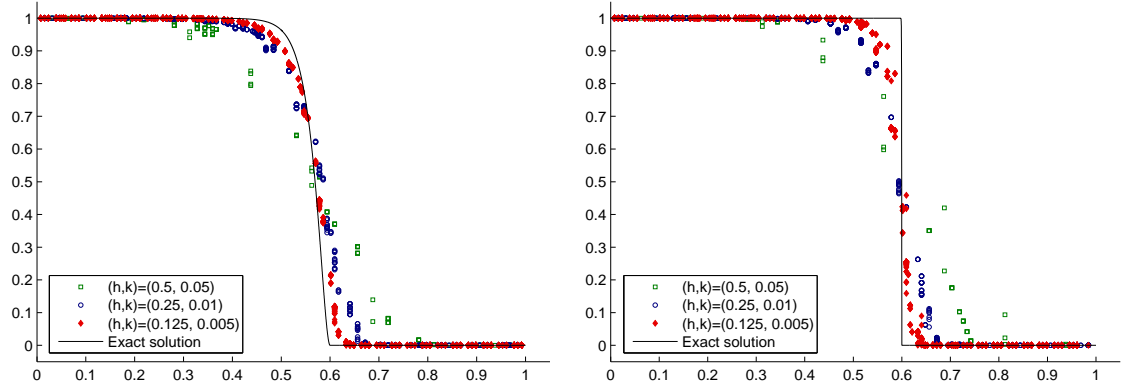}
\caption{The approximate solution profiles at the time $t=0.5$ for $\delta=0.01$ and $\delta=0.0001$}\label{fig:prof}
\end{figure}

\begin{Rq}\label{rq:upwind}
We remark that the nonlinear system (\ref{loc:test1}) (or a linear one arising during the Newton's procedure)
has a special structure; more specifically, for each $K\in\{1,\ldots,card(\M)\}$ the $K$-th equation does not contain any unknown
different from $\uK \ta (\us)_{\s\in\E_K}$ (here we denote by $K$ both the control volume and the index of the unknown $u_K$);
therefore one can algebraically eliminate $\uK$, so that the number of equations to solve becomes $card(\E)$.
\end{Rq}
Since we do not impose many constraints on the mesh (in particular it can be nonconforming), it
is not difficult to perform a local grid refinement. Finally note that
there is a possibility to reduce the number of unknowns by using a method
introduced in \cite{Eym_Gal_Her_2}; one can eliminate
the interior interface unknowns $(u_{\s})_{\s\in
\E_{int}}$ by expressing them as a consistent barycentric
combinations of the values $u_K$.


\end{document}